\documentclass[a4paper,fleqn]{cas-sc}

\usepackage[authoryear,longnamesfirst]{natbib}

\usepackage{graphicx}%
\usepackage{multirow}%
\usepackage{amsmath,amssymb,amsfonts}%
\usepackage{mathrsfs}%
\usepackage{appendix}
\usepackage[dvipsnames]{xcolor}
\usepackage{xcolor}%
\usepackage{textcomp}%
\usepackage{manyfoot}%
\usepackage{booktabs}%
\usepackage{algorithm}%
\usepackage{algorithmicx}%
\usepackage{algpseudocode}%
\usepackage{listings}%
\usepackage{cleveref}
\usepackage{subfig}

\DeclareMathOperator*{\argmin}{argmin}

\definecolor{orange}{RGB}{255,127,0}

\newcommand{\remove}[1]{{\textcolor{red}{}}}
\newcommand{\crossout}[1]{{\textcolor{red}{\st{}}}}

\newcommand{\bfb}{\mathbf b}

\newcommand{\bfw}{\mathbf w}

\newcommand{\bfbtrue}{\bfb^{\mathrm{true}}}

\newcommand{\bfs}{\mathbf s}
\newcommand{\bfq}{\mathbf q}

\newcommand{\bfr}{\mathbf r}

\newcommand{\bfu}{\mathbf u}
\newcommand{\bfv}{\mathbf v}
\newcommand{\bfx}{\mathbf x}

\newcommand{\bfy}{\mathbf y}

\newcommand{\bfh}{\mathbf h}

\newcommand{\bfe}{\mathbf e}

\newcommand{\bfxtrue}{\bfx^{\mathrm{true}}}

\newtheorem{example}{Example}%
\def\tsc#1{\csdef{#1}{\textsc{\lowercase{#1}}\xspace}}
\tsc{WGM}
\tsc{QE}
\tsc{EP}
\tsc{PMS}
\tsc{BEC}
\tsc{DE}

\begin{document}
\let\WriteBookmarks\relax
\def\floatpagepagefraction{1}
\def\textpagefraction{.001}
\shorttitle{}
\shortauthors{A. Alsubhi and R. A. Renaut}

\title [mode = title]{Hybrid GKB Methods for X-Ray Tomography Problems with Unmatched Back Projector}                      



\author[1]{Abdulmajeed Alsubhi}[orcid=0009-0002-1984-1304]
\cormark[1]
\ead{ahalsubhi@iu.edu.sa}


\affiliation[1]{organization={Department of Mathematics, Faculty of Science, Islamic University of Madinah},city={Madinah},
                postcode={42351}, 
                country={Saudi Arabia}}

\author[2]{Rosemary Renaut}[orcid=0000-0001-9296-0890]

\ead{renaut@asu.edu}


\affiliation[2]{organization={School of Mathematical and Statistical Sciences, Arizona State University},
                addressline={AZ}, 
                postcode={85287}, 
                city={Tempe},
                country={USA}}

\cortext[cor1]{Corresponding author}


\begin{abstract}
X-ray Computed Tomography (CT) is a widely used imaging modality in medical, industrial, and scientific applications. In practical CT implementations, forward, $A$, and back projection, $B$, operators are often constructed using different discretization schemes to improve computational efficiency on available software and hardware platforms. Consequently, the resulting projector pair is generally unmatched, meaning that the back projector is not the exact adjoint of the forward projector. This mismatch alters the mathematical properties of the reconstruction problem and can affect the convergence behavior of iterative reconstruction methods.
Previous studies have proposed AB- and BA-Golub–Kahan bidiagonalization (GKB) methods, as well as GMRES and hybrid GMRES methods, for solving CT reconstruction problems with unmatched projector pairs, demonstrating reduced semiconvergence effects compared with conventional approaches. In this work, we develop hybrid variants of the AB- and BA-GKB algorithms by incorporating regularization within the projection process. The singular value decomposition for the operator on the projected space is used to efficiently reconstruct the solution, and to automatically select the regularization parameter using either the L-curve or generalized cross-validation.  Numerical experiments on several CT reconstruction problems demonstrate the effectiveness of the proposed hybrid methods in improving robustness against semiconvergence.
\end{abstract}




\maketitle

\section{Introduction}\label{sec:intro}
X-ray computed tomography (CT) reconstructs an image of an object from X-ray measurements acquired at multiple projection angles. The reconstruction problem can be formulated as
\begin{equation}\label{eq:inv prob}
    A\bfx \approx \bfb.
\end{equation}
Here $\bfx \in \mathbb{R}^{n}$ is the unknown image, $\bfb \in \mathbb{R}^{m}$ contains the measured projection data assumed to be contaminated by additive white Gaussian noise, and 
$A \in \mathbb{R}^{m\times n}$ denotes the forward projection operator with no restrictions on its dimensions or rank.

The least squares problem for \eqref{eq:inv prob}
\begin{equation*}\label{eq:LS}
    \min_{\bfx \in \mathbb{R}^n}\|A\bfx-\bfb \|_2,
\end{equation*}
is solved by the normal equations,   either
\begin{equation} \label{eq:over NE}
   A^\top A\bfx =A^\top \bfb,
\end{equation}
for  $m\geq n$, or 
\begin{equation} \label{eq:under NE}
   A A^\top \bfy = \bfb, \ \ \ \bfx=A^\top \bfy,
\end{equation}
for the underdetermined case, $m<n$. Both situations can be solved using standard iterative approaches,  but for practical CT applications  the dimensions of $A$ are typically very large, making its explicit construction and storage computationally prohibitive.  
Consequently, iterative reconstruction algorithms are commonly implemented in a matrix-free framework, in which only matrix-vector products with $A$ and its transpose are required. These operations are performed using numerical routines that approximate the underlying imaging physics via appropriate discretization schemes and the forward and back projection operators are often constructed using different discretizations and software implementations. As a result, the back projector $B \in \mathbb{R}^{n\times m}$ is generally not the exact adjoint $A^\top$ of the forward projector, that is, $B \neq A^{\top}$. Then, the pair $A,B$ is referred to as an unmatched projector pair \citep{article,hansen2022gmres} and the normal equations are replaced with the generally unsymmetric unmatched normal equations, 
\begin{equation} \label{eq:over UN NE}
   B A\bfx =B \bfb,
\end{equation}
and
\begin{equation} \label{eq:under UN NE}
   A B \bfy = \bfb, \ \ \ \bfx=B \bfy,
\end{equation}
respectively, for $m \ge n$ and $m<n$. If the range of $A^\top$ equals the range of $B$, then the solutions of the matched and unmatched equations are equivalent \citet{wathen2025least}. Here, we consider unmatched pairs that do not satisfy this range condition. 

The use of unmatched projector pairs fundamentally alters the structure of the normal equations in \cref{eq:over NE,eq:under NE} and classical iterative solvers, such as the Conjugate Gradient algorithm, do not necessarily yield reliable reconstructions  \citep{article,doi:10.1137/17M1133828}. The AB- and BA-GMRES methods were introduced in \citet{doi:10.1137/070696313} to solve general inverse problems with nonsymmetric operators, and several choices of the operator $B$, including diagonal scaling strategies were designed to improve convergence. The AB- and BA-GMRES were adapted for CT reconstruction problems with unmatched projector pairs in \citet{hansen2022gmres}. There, GMRES is applied to the unmatched normal equations without requiring the explicit construction of the adjoint operator. Often exhibiting rapid initial convergence, these approaches remain susceptible to semiconvergence, in which early iterates approach a desirable approximation of the true solution before eventually becoming contaminated by propagated noise.

Incorporating Tikhonov regularization within the Krylov subspace updates  reduces the issues of semiconvergence, providing the hybrid variants of AB- and BA-GMRES,  \citet{bentley2026hybrid}.  Their effectiveness depends strongly on the properties of the projected problem generated by the Arnoldi process for which the Hessenberg matrix for the projected system governs both the regularization procedure and the quality of the reconstructed solution \citet{saad2003iterative}.

An alternative Krylov subspace framework is provided by the GKB process, originally developed by Golub and Kahan \citet{golub1965calculating} for reducing a general matrix to bidiagonal form. Besides its importance in computing the singular value decomposition (SVD), GKB serves as the foundation of several widely used iterative methods, including LSQR \citet{paige1982lsqr} and LSMR \citet{fong2011lsmr}. Additional applications of GKB-based techniques can be found in \citep{bjorck1988bidiagonalization,bentbib2018solution,beik2020golub,bianchi2026iterated}. Moreover, AB- and BA-GKB methods for reconstructing CT problems were presented in \citet{alsubhi2026GKB}. Unlike GMRES, which generates Hessenberg matrices, GKB produces lower bidiagonal matrices that are considerably sparser and retain important spectral information from the original problem. As a consequence, AB- and BA-GKB have been observed to exhibit milder semiconvergence and improved reconstruction quality compared with their GMRES counterparts \citet{alsubhi2026GKB}.

Motivated by these observations, we develop hybrid AB- and BA-GKB methods for CT reconstruction with unmatched projector pairs. The bidiagonal structure of the projected space operator generated by GKB, assumed to be small, enables efficient construction of its SVD and facilitates efficient regularization parameter determination using standard methods. As examples,  we consider L-curve and GCV for parameter estimation. Since regularization is applied only to a projected problem of relatively small dimension, the computational cost associated with the SVD and parameter-selection procedure is negligible compared with the matrix-vector operations required by the GKB iterations.  Numerical results demonstrate that the proposed hybrid AB- and BA-GKB methods significantly reduce the effects of semiconvergence and can yield more accurate reconstructions than their GMRES-based counterparts.

\subsection{Major Contributions}\label{sec:major}
\begin{itemize}
    \item Introduction of hybrid AB-GKB and BA-GKB algorithms for CT reconstruction with unmatched forward and backprojectors.
    \item Stabilization of the projected space operator using standard Tikhonov regularization,  with GCV- and L-curve-based parameter selection using the SVD of the projected space operator.
     \item Numerical validation demonstrating that the proposed methods achieve reconstruction quality comparable to that of state-of-the-art hybrid GMRES approaches.
\end{itemize}

\subsection{Organization }\label{sec:organ}
The organization of this paper is outlined as follows. \Cref{sec:Background} reviews the mathematical techniques used to solve inverse problems with unmatched pairs using GMRES, covering both classical and hybrid approaches. \Cref{sec:ABBAGKB} revisits the GKB algorithm and introduces  new $AB$ and $BA$ hybrid methodologies. \Cref{sec:Numerical examples} presents our numerical examples, while \Cref{sec:Conclusion} concludes our study and highlights potential directions for future research.

\section{AB- and BA- Iterative Methods for Unmatched Normal Equations} \label{sec:Background}
AB- and BA-GMRES methods for unmatched pairs were introduced in \citet{hansen2022gmres} and extended in the form as hybrid methods to include iterative regularization by \citet{bentley2026hybrid}. Here we briefly review these algorithms for comparison with an alternative formulation and extension based on the GKB algorithm, \citet{alsubhi2026GKB}. 

AB-GMRES and BA-GMRES find at iteration $k$ the least squares solution of $\|AB\bfy-\bfb \|^2_2$, with update $\bfx=B\bfy$, and of $\|BA\bfx-\bfb \|^2_2$ for the unmatched normal equations \cref{eq:over UN NE,eq:under NE}, respectively. Both algorithms use the same Krylov space $\mathcal{K}_k(BA, B\bfb)$ and construct an orthogonal basis for the approximate solution with least squares residual using an Arnoldi iteration.   Let $M$ be the coefficient matrix $AB$ or $BA$, then $k$ applications of Arnoldi iteration yield
\begin{equation*}\label{eq:Arn iter}
    M \bar{Q}_k=\bar{Q}_{k+1}H_k,
\end{equation*}
where $\bar{Q}_{k+1}= [\bar{\bfq}_1,\bar{\bfq}_2,\dots,\bar{\bfq}_{k+1}]\ $ and the matrix $H_{k+1,k}$ has upper Hessenberg structure, 
\begin{equation*}
     H_k=
\begin{bmatrix}\label{upper Hessenberg}
   h_{1,1}&h_{1,2} &\dots & h_{1,k}\\
   h_{2,1}&h_{2,2} &\dots & h_{2,k}\\
   &\ddots &\ddots &\vdots\\
   &&h_{k,k-1}&h_{k,k}\\
   &&&h_{k+1,k}
\end{bmatrix}.
\end{equation*}
Using this decomposition in  \cref{eq:over UN NE,eq:under UN NE} yields the projected problem 
\begin{equation}\label{eq:proj}
    \bfu_k=\argmin_{\bfu} \Big\| H_k \bfu - \bfe_1 \|\bfr_0\|_2 \Big\|_2,
\end{equation}
where $\bfe_1 =(1,0,\dots,0)^\top \in \mathbb{R}^{k+1}$.

For ill-conditioned problems, for which the underlying singular spectrum is broad, Krylov methods demonstrate semi-convergence. The initial iterations capture the dominant singular spectrum and the solution converges to the exact solution. But as the iteration proceeds, the Krylov subspace starts to include singular vectors, associated with smaller singular values, which are often contaminated by noise present in the data, and the solution departs from the exact solution. This problem is exacerbated for unmatched $AB$ and $BA$ operators for which the spectrum not only does not match the matched operator spectrum but is also unstable due to non-normality. Therefore, standard iterative methods applied to \eqref{eq:proj} typically require appropriate stopping rules to stop the iterations before the onset of semi-convergence. Then the number of iterations serves as a regularization parameter, \citet{hansen2021stopping}. Alternatively, one may seek a regularized solution on the projected space within the iteration, \citet{ChungGazzolareview}. Applied to the AB- and BA-GMRES algorithms  Tikhonov regularization on the projected space \eqref{eq:proj} yields the $\lambda$-dependent solution 
\begin{equation}\label{eq:ABGMRES Hyb}
    \bfu_{k,\lambda}=\argmin_{\bfu}\left \|\begin{bmatrix}
        H_k\\\lambda I
    \end{bmatrix} \bfu - \|\bfr_0\|_2\begin{bmatrix}
        \bfe_1\\ 0
    \end{bmatrix}  \right\|^2_2. 
\end{equation}
An equivalent solution to \eqref{eq:ABGMRES Hyb} can be given using the solution of the normal equations
\begin{equation}\label{eq:normal ABGMRES Hyb}
     \bfu_{k,\lambda}= \|\bfr_0\|_2(H^\top_k H_k +\lambda^2I)^{-1} H^\top_k
        \bfe_1  .
\end{equation}
An option for the solution of \eqref{eq:normal ABGMRES Hyb} is to calculate the singular value decomposition of $H_k$ and then use a standard parameter choice method to find $\lambda$, \citet{ChungGazzolareview}. 

The hybrid AB-GMRES and hybrid BA-GMRES Algorithms as detailed in \citet{bentley2026hybrid} are summarized in \Cref{alg:ABGMRES,alg:BAGMRES}.

\begin{minipage}{0.47\textwidth}
\begin{algorithm}[H]
    \centering
    \caption{Hybrid AB-GMRES \citet{doi:10.1137/070696313}}\label{alg:ABGMRES}
    \begin{algorithmic}[1]
       \Require $A$, $B$, $\bfb$, 
       \State \text{Choose initial $\bfy_0$}, \text{set ${\bfx}_0=B\bfy_0$}
        \State {$\bfr_0=\bfb-A\bfx_0$}
        \State {$\bar{\bfq}_1=\bfr_0 / \| \bfr_0\|_2$}
        \For{$k=1,2,\dots$}
        \State ${\bar{\bfw}_k=A(B \bar{\bfq}_k})$ 
        \For{$i=1,2,\dots,k$}
        \State ${\bfh_{i,k} = \bar{\bfq}^\top_i \bar{\bfw}_k}$ 
        \State ${\bar{\bfw}_k= \bar{\bfw}_k- \bfh_{i,k} \bar{\bfq}_i }$ 
         \EndFor
        \State ${\bfh_{k+1,k} = \|\bar{\bfw}_k \|_2}$ 
        \State ${\bar{\bfq}_{k+1}=\bar{\bfw}_k / \bfh_{k+1,k}}$ 
        \State Determine $\lambda$
        \State {Update $\bfu_{k,\lambda}$: solve \eqref{eq:normal ABGMRES Hyb}}
        \State {$\bfx_k=\bfx_0+B( \bar{Q}_k \bfu_k$)}
        \State \text{$\bfr_k=\bfb-A\bfx_k$}
        \State \text{Stopping rule is applied here}
        \EndFor
    \end{algorithmic}
\end{algorithm}
\end{minipage}
\hfill
\begin{minipage}{0.47\textwidth}
\begin{algorithm}[H]
    \centering
    \caption{Hybrid BA-GMRES \citet{doi:10.1137/070696313}}\label{alg:BAGMRES}
    \begin{algorithmic}[1]
      \Require $A$; $B$, $\bfb$, 
        \State \text{Choose initial $\bfx_0$} 
        \State {$\bfr_0=B(\bfb-A\bfx_0)$}
        \State {$\bar{\bfq}_1=\bfr_0 / \| \bfr_0\|_2$}
        \For{$k=1,2,\dots$}
        \State ${\bar{\bfw}_k=B (A \bar{\bfq}_k})$ 
        \For{$i=1,2,\dots,k$}
        \State ${\bfh_{i,k} = \bar{\bfq}^\top_i \bar{\bfw}_k}$ 
        \State ${\bar{\bfw}_k= \bar{\bfw}_k- \bfh_{i,k} \bar{\bfq}_i }$ 
         \EndFor
        \State ${\bfh_{k+1,k} = \|\bar{\bfw}_k \|_2}$ 
        \State ${\bar{\bfq}_{k+1}=\bar{\bfw}_k / \bfh_{k+1,k}}$ 
        \State Determine $\lambda$
        \State {Update $\bfu_{k,\lambda}$: solve \eqref{eq:normal ABGMRES Hyb}}
        \State {$\bfx_k=\bfx_0+ \bar{Q}_k \bfu_k$}
        \State \text{$\bfr_k=B(\bfb-A\bfx_k)$}
        \State \text{Stopping rule is applied here}
        \EndFor
    \end{algorithmic}
\end{algorithm}
\end{minipage}

The system matrix, $H_k$, in  \Cref{alg:ABGMRES,alg:BAGMRES} has potentially $(k^2+3k)/2$  non-zero elements and  as $k$ increases the computation of the SVD may become  computationally expensive. In the next section, we review and extend an approach based on the use of the Golub Kahan Bidiagonalization (GKB) in place of GMRES, as introduced in \citet{alsubhi2026GKB}.

\section{Hybrid AB- and BA-GKB Algorithms}\label{sec:ABBAGKB}
Following the work on the use of AB- and BA-GMRES for unmatched operators, AB- and BA-GKB were developed in \citet{alsubhi2026GKB}. Underlying the GKB algorithm is the decomposition   
\begin{equation*}\label{eq:GKB dec}
    A Q_k=S_{k+1} C_k,
\end{equation*}
where $Q_{k}$ and $S_{k+1}$ are column orthogonal matrices and
\begin{equation*}\label{eq:lower bid}
   C_{k} = \begin{bmatrix}
       \alpha_1&&&&\\
   t_2&\alpha_2&&&\\
   &\ddots&\ddots&&\\
      &&&t_k&\alpha_k\\
     &&&&t_{k+1}
   \end{bmatrix}
\end{equation*} 
is a lower bidiagonal matrix. 
For unmatched problems \cref{eq:over UN NE,eq:under UN NE}, the GKB yields the projected problem, 
$\min_{\bfu \in \mathbb{R}^k}\|C_k\bfu- t_1 \bfe_1\|_2^2$, 
corresponding to direct application of the GKB algorithm to the system matrix $A$ with right, or left preconditioner $B$, respectively. The solution is obtained as $\bfx = Q_k\bfu_k$ in each case, but the AB- algorithm is initialized with $\bfb/\|\bfb\|_2$ and the BA- algorithm with $B\bfb/\|B\bfb\|_2$. 

As for the algorithms based on the GMRES, an approach using the GKB is subject to semiconvergence of the iterations for ill-conditioned problems. Under the same framework with appropriate stopping conditions applied to GMRES and GKB algorithms for unmatched operators, \citet{alsubhi2026GKB} showed that the AB- and BA-GKB methods outperform their GMRES-based counterparts and are less susceptible to semiconvergence. To further improve reconstruction quality, we incorporate a regularization term into the AB- and BA-GKB methods yielding their hybrid counterparts based on the iterative solution of 
\begin{equation}
  \bfu_{k,\lambda}=\argmin_{\bfu}  \left\|\begin{bmatrix} \label{eq:ABGKB Hyb}
            C_k \\ \lambda I
        \end{bmatrix}  \bfu - t_1 \begin{bmatrix}
            \bfe_1\\ 0
        \end{bmatrix} \right\|^2_2,
\end{equation}
equivalent to \eqref{eq:ABGMRES Hyb} but with bidiagonal $C_k$ replacing upper Hessenberg $H_k$, and where $t_1=\|\bfb\|_2$ for the AB methods and $\|B\bfb\|_2$ for the BA methods. The equivalent regularized solution to \eqref{eq:ABGKB Hyb} based on the normal equations can be written as
\begin{equation}\label{eq:normal ABGKB Hyb}
     \bfu_{k,\lambda}= t_1(C^\top_k C_k +\lambda^2I)^{-1} C^\top_k 
        \bfe_1  .
\end{equation}

The hybrid AB-GKB and BA-GKB algorithms  are summarized in \Cref{alg: Hybrid AB-GKB,alg: Hybrid BA-GKB}. 

\begin{minipage}{0.47\textwidth}
\begin{algorithm}[H]
    \centering
    \caption{Hybrid AB-GKB}\label{alg: Hybrid AB-GKB}
    \begin{algorithmic}[1]
 \Require $A$, $B$,  $\bfb$
       \State $t_1=\|\bfb\|_2$ 
       \State $\bfs_1=\bfb/t_1$, \ $\bfq=B^\top (A^\top \bfs_1)$
       \State $\alpha_1=\|\bfq\|_2$
       \State $\bfq_1=\bfq/\alpha_1$
        \For{$k=1,2,\dots$}
\State {$\bfs_k=A(B\bfq_{k-1})- \alpha_{k-1} \bfs_{k-1}$ }
\State {$t_k=\|\bfs_k \|_2$, \ $\bfs_k=\bfs_k /t_k $ }
\State {$\bfq_k=B^\top (A^\top \bfs_k)- t_k \bfq_{k-1}$ \;}
\State {$\bfq_k = \bfq_k - \sum^{k-1}_{i=1} (\bfq^\top_i \bfq_k)  \bfq_i$}
\State {$\alpha_k=\|\bfq_k \|_2, \ \bfq_k=\bfq_k/\alpha_k$ \;}
\State{$\bfs_i$ and $\bfq_i$ define $S_{k+1}$,  $Q_k$}
\State{Scalars $t_k$ and $\alpha_k$ define  $C_k$}       
        \State \text{Compute} $\lambda$
         \State {Update $\bfu_{k,\lambda}$: solve \eqref{eq:normal ABGKB Hyb}}
        \State {$\bfx_k=B (Q_k \bfu_k$)}
        \State \text{Stopping rule is applied here}
        \EndFor
    \end{algorithmic}
\end{algorithm}
\end{minipage}
\hfill
\begin{minipage}{0.47\textwidth}
\begin{algorithm}[H]
    \centering
    \caption{Hybrid BA-GKB}\label{alg: Hybrid BA-GKB}
    \begin{algorithmic}[1]
                \Require $A$, $B$,  $\bfb$
       \State $t_1=\|B\bfb\|_2$ 
       \State $\bfs_1=B\bfb/t_1$, \ $\bfq=A^\top (B^\top  \bfs_1)$
       \State $\alpha_1=\|\bfq\|_2$
       \State $\bfq_1=\bfq/\alpha_1$
        \For{$k=1,2,\dots$}
\State {$\bfs_k=B(A\bfq_{k-1})- \alpha_{k-1} \bfs_{k-1}$ }
\State {$t_k=\|\bfs_k \|_2$, \ $\bfs_k=\bfs_k /t_k $ }
\State {$\bfq_k=A^\top (B^\top \bfs_k)- t_k \bfq_{k-1}$ \;}
\State {$\bfq_k = \bfq_k - \sum^{k-1}_{i=1} (\bfq^\top_i \bfq_k)  \bfq_i$}
\State {$\alpha_k=\|\bfq_k \|_2, \ \bfq_k=\bfq_k/\alpha_k$ \;}
\State{$\bfs_i$ and $\bfq_i$ define $S_{k+1}$,  $Q_k$}
\State{Scalars $t_k$ and $\alpha_k$ define  $C_k$}       
        \State \text{Compute} $\lambda$
       \State {Update $\bfu_{k,\lambda}$: solve \eqref{eq:normal ABGKB Hyb}}
        \State {$\bfx_k=Q_k \bfu_k$}
        \State \text{Stopping rule is applied here}
        \EndFor
    \end{algorithmic}
\end{algorithm}
\end{minipage}

\subsection{Computational Costs of GMRES and GKB Algorithms} \label{sec:comp costs}
The computational costs of the non-hybrid GMRES and GKB algorithms were presented in \citet{alsubhi2026GKB} and are dominated by the costs of the matrix-vector operations.  In particular,  both the AB-GMRES and the AB-GKB algorithms require an extra matrix vector multiplication in the update of $\bfx_k$, steps 14 and 15 in \Cref{alg:ABGMRES,alg: Hybrid AB-GKB}, respectively.  Thus, as compared to the respective BA algorithms the computational complexity increases by $2mnk$. Supposing in each case, that the SVD of $H_k$, respectively $C_k$, is used to find the solution $\bfu_{k,\lambda}$, then the difference between the GKB and GMRES algorithms is impacted by the structure of these matrices, the SVD for $H_k$ being more expensive. As indicated in \citet{alsubhi2026GKB} the computational costs obtained are 

\begin{equation}\label{cost:ABGMRES}
 \text{Cost with AB-GMRES :}  \mathcal{O}(6mnk+4mk^2+[(k^4+3k^2)/2]),
\end{equation}
 
\begin{equation}\label{cost:BAGMRES}
  \text{Cost with BA-GMRES  :}  \mathcal{O}(4mnk+4nk^2 + [(k^4+3k^2)/2] ),
\end{equation}

\begin{equation}\label{cost:ABGKB}
  \text{Cost with AB-GKB :}  \mathcal{O}(10mnk+4mk^2+[2k^3]),
\end{equation}

\begin{equation}\label{cost:BAGKB}
  \text{Cost with  BA-GKB :}  \mathcal{O}(8mnk+4nk^2+ [2k^3]).
\end{equation}
While the terms in parentheses \cref{cost:ABGMRES,cost:BAGMRES} and \cref{cost:ABGKB,cost:BAGKB}, for the calculation of the SVD for $H_k$ and $C_k$, respectively,  are small compared to large $m$ and $n$ for a limited number of iterations $k$, it is evident that the SVD cost of $C_k$ is approximately a factor of $1/k$ smaller than that of $H_k$.

Assuming that the algorithms stop with $k<<[m,n]$, \crefrange{cost:ABGMRES}{cost:BAGKB} show that the GKB-based algorithms are respectively, five thirds and two times more computationally expensive that their GMRES counterparts for a fixed choice of $k$. This was confirmed in \citet{alsubhi2026GKB} for the non-hybrid algorithms, where it was shown that the number of iterations to satisfaction of equivalent stopping rules was higher for the GKB algorithms and that their cost was higher. It is therefore of interest to investigate the performance of the comparative hybrid algorithms for the inversion of large scale problems. First, this requires a description of the approach used to find the regularization parameter $\lambda$ at each iteration of each algorithm.

\subsection{Selection of Regularization Parameters}
The selection of $\lambda$ in \Crefrange{alg:ABGMRES}{alg: Hybrid BA-GKB} is important for the quality of the reconstruction. If the noise variance is known \textit{a priori}, options to determine $\lambda$ include the discrepancy principle or the unbiased predictive risk estimator. When the noise variance is unknown, methods for estimating $\lambda$ include the L-curve and GCV. For consistency with real applications that often lack sufficient information about the noise variance, we consider the GCV, \citet{golub1979generalized}, and the L-curve, \citet{hansen1993use}. The motivation behind these choices is that our numerical experiments show that these methods effectively mitigate the impact of noise on reconstructions. 

In our experiments we assume that $k$ is small enough that we can use the SVD of $C_k$, respectively $H_k$, to evaluate the parameter-estimation techniques. We present the relevant formulae for the system with matrix $C_k$ and right hand side $t_1 \bfe_1$, noting that the formulae apply for $H_k$ with right hand side $\|\bfr_0\|_2\bfe_1$ in the same way. Now, let $C_k = U \Sigma V$ denote the full SVD of $C_k$, where $U =[\bfu_1, \bfu_2, \dots,\bfu_k]$ and $V=[\bfv_1, \bfv_2, \dots,\bfv_k]$ are orthogonal matrices and $\Sigma$ is a (possibly rectangular) diagonal matrix containing the singular values $\sigma_j$ on its main diagonal in descending order. The SVD filter factors take the form $f_i(\lambda) = \sigma^2_i/(\sigma^2_i+\lambda^2)$, for $i=1,2,\dots,k$, and  
the regularized solution $\bfu_{k,\lambda}$, given in \eqref{eq:normal ABGKB Hyb}, can be written as 
\begin{equation*}\label{eq:solnSVD}\bfu_{k,\lambda} = \sum^{k}_{i=1}t_1 f_i(\lambda) \frac{\bfu^\top_i\bfe_1}{\sigma_i}\bfv_i.
\end{equation*}

The GCV determines the minimum value of the GCV function, here given immediately in terms of the SVD as 
\begin{equation}\label{eq:GCV SVD}
    GCV(\lambda) = \frac{\sum^{k}_{i=1} t_1^2f^2_i(\lambda)  (\bfu_i^\top\bfe_1)^2}{\sum^{k}_{i=1} f^2_i(\lambda)}.
\end{equation}
Even though the GCV is effective at finding a good estimate for $\lambda$ for most problems, it has some difficulty when the GCV function is flat, \citet{varah1983pitfalls}. This is because the noise in the data is small and the numerator in \eqref{eq:GCV SVD} is approximately constant for a large range of $\lambda$. We consider, however, a reasonable amount of noise and the scenario of a flat GCV function is not a concern.

We also consider the L-curve in which $\lambda$ is selected as the corner of the curve generated by plotting the log values of the following two quantities against each other, again in terms of the SVD 
 \begin{equation}\label{eq:SVD Lcurve}
   \left( \sqrt{\sum^{k}_{i=1} t_1^2f_i(\lambda)^2  (\bfu^\top_i\bfe_1)^2} , \sqrt{\sum^{k}_{i=1} f_i(\lambda) t_1\frac{\bfu^\top_i\bfe_1}{\sigma_i}\bfv_i} \ \right).
\end{equation}
The corner in the L-curve indicates that the solution $\bfu_k$ reaches a suitable balance between the noise remaining in $\|\bfu_{k,\lambda}\|_2$ and the residual norm $\|C_k\bfu_{k,\lambda}- t_1 \bfe_1\|_2$. 

\Cref{eq:GCV SVD,eq:SVD Lcurve} facilitate efficient calculation of the regularization parameter over a range of values, based on the approximate singular values of the projected system matrix.

\section{Numerical Experiments}\label{sec:Numerical examples}
In this section, we assess the performance of the proposed hybrid AB-GKB and BA-GKB algorithms for reconstructing CT problems with unmatched forward and backward projector pairs. We compare the reconstructions produced by these algorithms with those derived from hybrid AB-GMRES and BA-GMRES. For all simulations, we assume the availability of both the exact image $\bfxtrue$ and the noise-free data $\bfbtrue$. To evaluate the impact of noise, we introduce additive white Gaussian noise to $\bfbtrue$. The Signal-to-Noise Ratio (SNR) is employed to quantify the noise present in $\bfb$.

To evaluate the quality of the reconstruction, we record the relative reconstruction error (RRE) between the exact image $\bfxtrue$ and the reconstructed image $\bfx^{(k)}$ at iteration $k$, as given by   
\begin{equation}\label{eq:RE}
    \text{RRE}(\bfx_{k})=\frac{\|\bfx_{k}- \bfxtrue\|_2} {\|\bfxtrue\|_2}.
\end{equation}
In most practical problems, the exact image is not available; hence, using the RRE to stop the iterations at the point of semiconvergence is not possible. Here, we assume an accurate estimation of the noise level $\|\bfe\|_2=\|\bfbtrue-\bfb\|_2$ is available and we terminate the iteration using the DP, \citet{morozov2012methods}, as 
\begin{equation}\label{eq:DP}
    k_{DP}=\min_k{\|A\bfx_k-\bfb\|_2\leq \tau \|\bfe\|_2},
\end{equation}
where $\tau\geq 1$ is a safety factor set to $1$.

The Residual Norm Stagnation (RNS) occurs when the solver stops making meaningful progress. The iteration is terminated using RNS when 
\begin{equation}\label{eq:RNS}
    k_{RNS}= \frac{ \bigr| \|A\bfx_{k-1}-\bfb\|_2 - \|A\bfx_k-\bfb\|_2 \bigr|}  {\|A\bfx_{k-1}-\bfb\|_2} < \epsilon,
\end{equation}
where the tolerance $\epsilon$ is user-defined to ensure convergence, and its appropriate choice may vary across algorithms. To facilitate a consistent comparison between the proposed methods and GMRES-based approaches, $\epsilon$ is selected empirically for each simulation—based on numerical experiments (not reported)—to minimize the RRE for both hybrid AB-GMRES and BA-GMRES algorithms.

For each problem, we compare the results obtained using the DP in \eqref{eq:DP} with those from the  RNS in \eqref{eq:RNS}, with the RNS results always given in parentheses in the tables. In our implementations, we record the RRE and the corresponding iteration counts. Furthermore, we report average timings across $5$ runs for $150$ iterations to validate the theoretical computational discussion presented in \Cref{sec:comp costs}.  For each problem we also contrast the results using the GCV and L-curve for determination of the regularization parameter.

All forward and backward projectors used in this section are generated with the GPU-accelerated ASTRA Toolbox \citet{van2015astra}. The operators $A$ and $B$ are implemented in a matrix-free manner, avoiding the explicit construction and storage of the corresponding matrices. Furthermore, the ASTRA Toolbox is extended to support the transpose operations $A^\top$ and $B^\top$, which are required for implementing the AB-GKB and BA-GKB algorithms.

All numerical results presented were computed using MATLAB Version 2026a \citet{MATLAB} on a desktop equipped with an Intel(R) Core(TM) i5-14400F CPU, 16 GB of memory, and an NVIDIA GeForce RTX 5060 GPU. The software will be available on request to the author.

\begin{figure}
  \centering
   \subfloat[True image]{\label{Fig:True mri}\includegraphics[width=.20\textwidth]{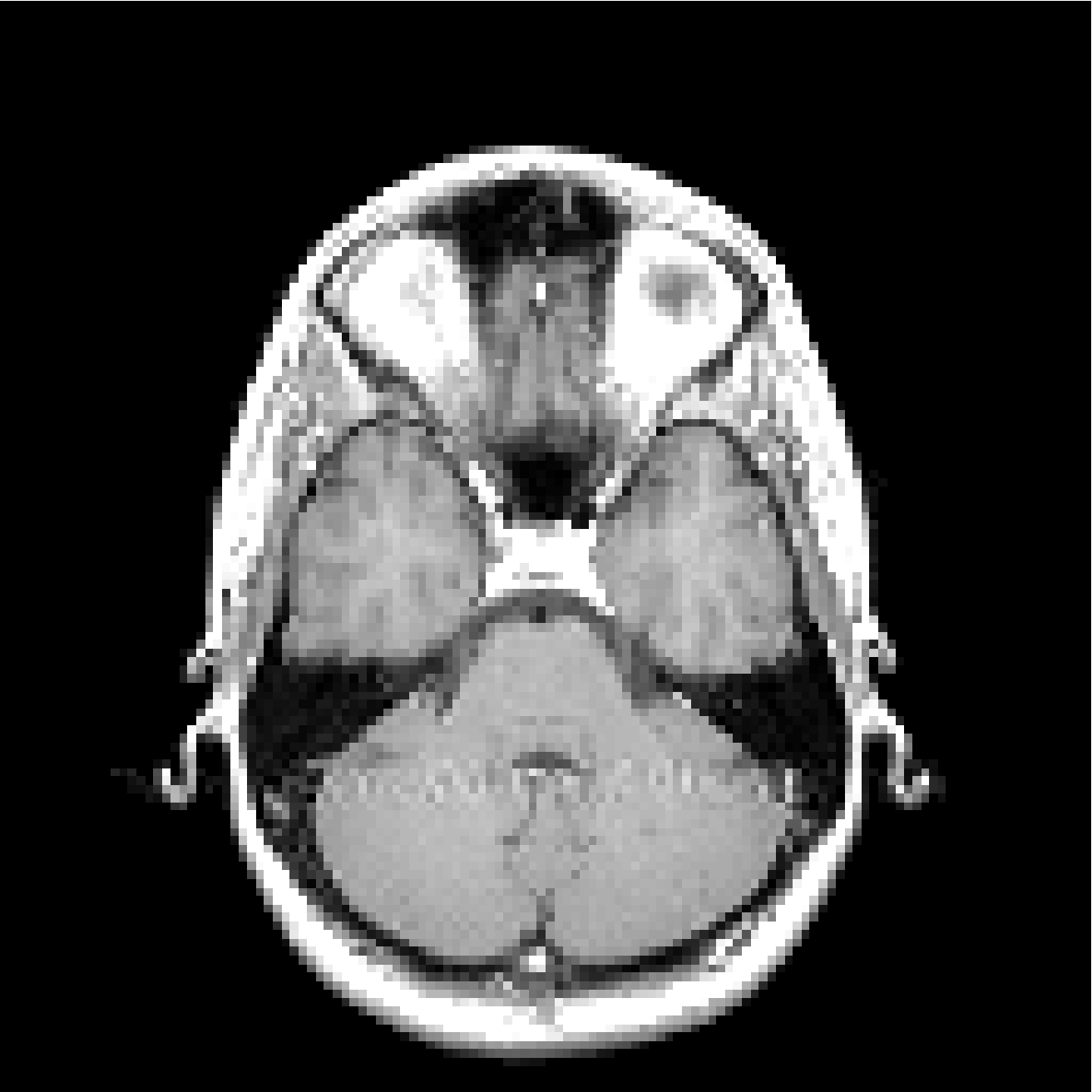}} \
   \subfloat[Noisy sinogram]{\label{Fig:blurred mri}\includegraphics[width=.20\textwidth]{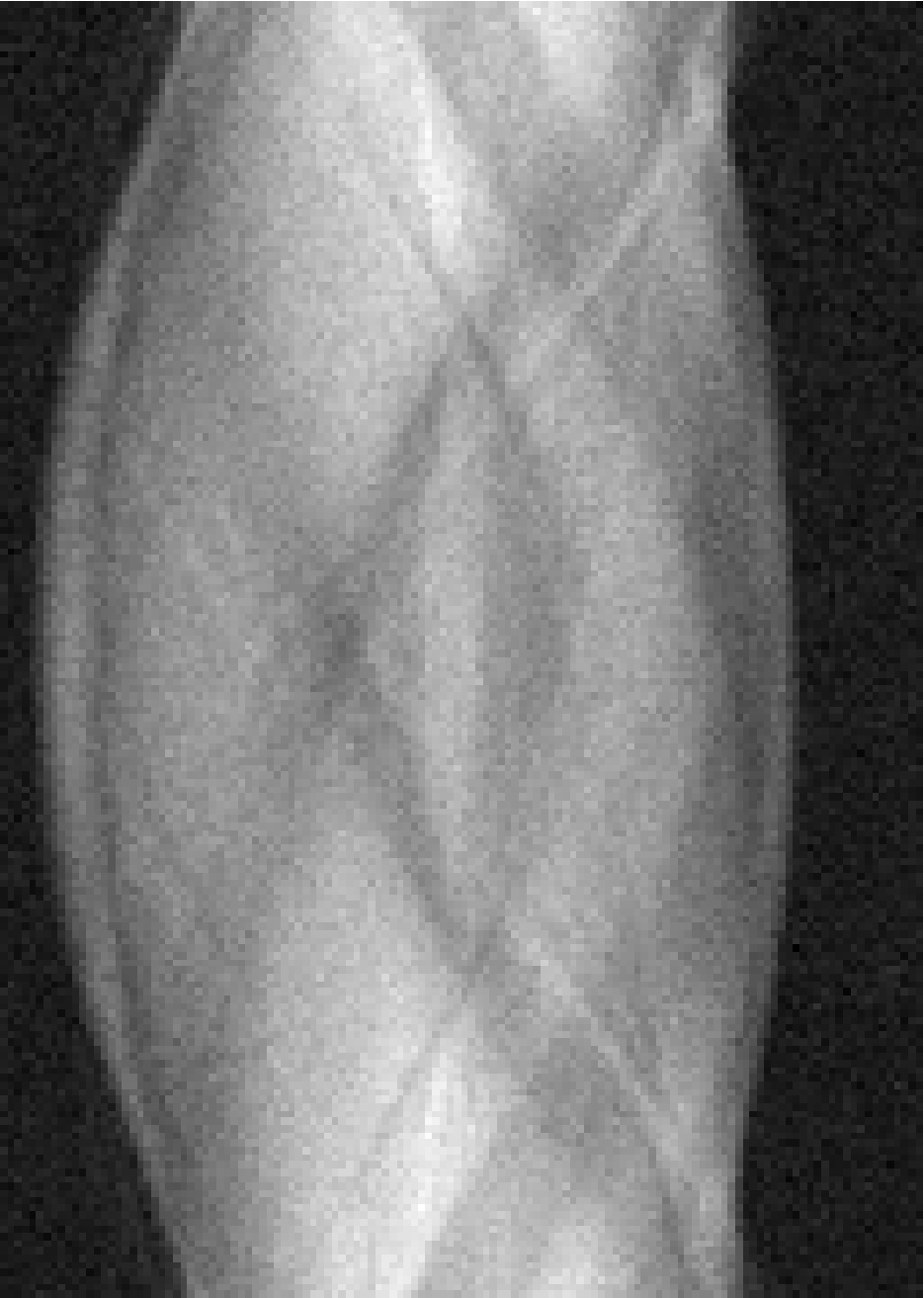}}
 \caption{The true image of $128 \times 128$ pixels and noisy sinogram with $6\%$ Gaussian noise, for the problem in \Cref{ex:problem1}. \label{Fig:mri True and noisy}}  
  \end{figure}

\begin{figure}
 \centering
  \subfloat[GCV]{\label{Fig:mri RRE GCV}\includegraphics[width=.48\textwidth]{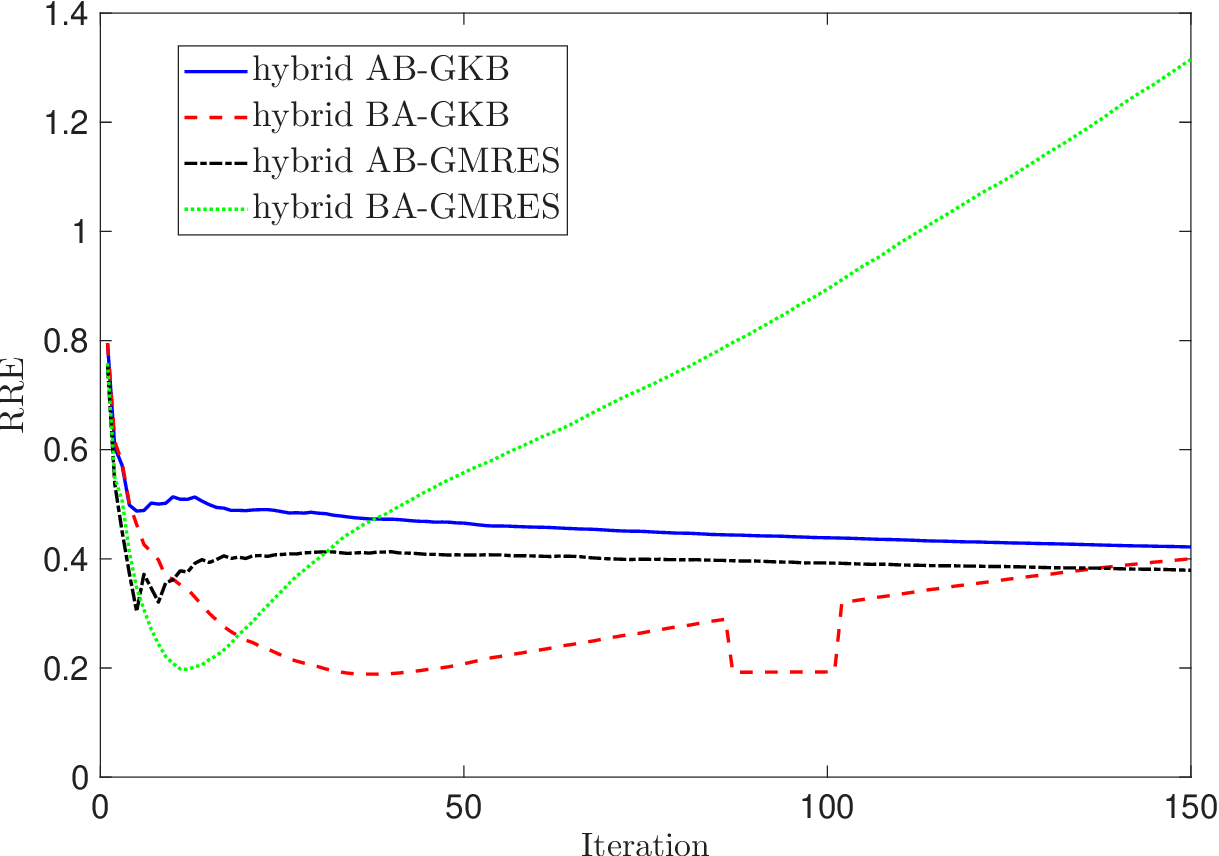}}   \ 
 \subfloat[L-curve]{\label{Fig:mri RRE Lcurve}\includegraphics[width=.48\textwidth]{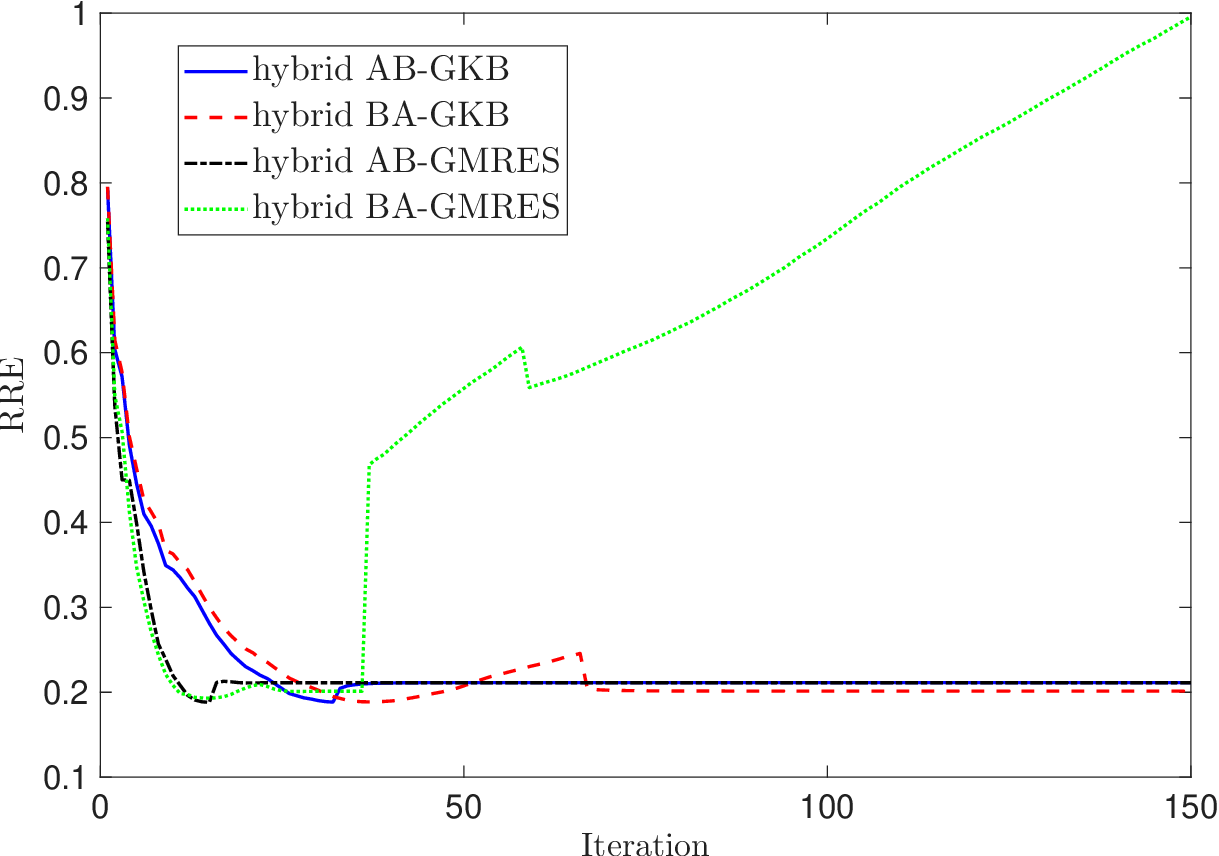}} 
 \caption{RRE for different methods for unmatched projectors  in \Cref{ex:problem1}.}\label{Fig:Small}  
 \end{figure}

\begin{figure}
  \centering
   \begin{tabular}{cccc} 
     \textbf{GCV, DP} & \textbf{GCV, RNS}  & \textbf{L-curve, DP} & \textbf{L-curve, RNS}  \\
    {\label{Fig:mri ABGKB GCV DP}\includegraphics[width=.18\textwidth]{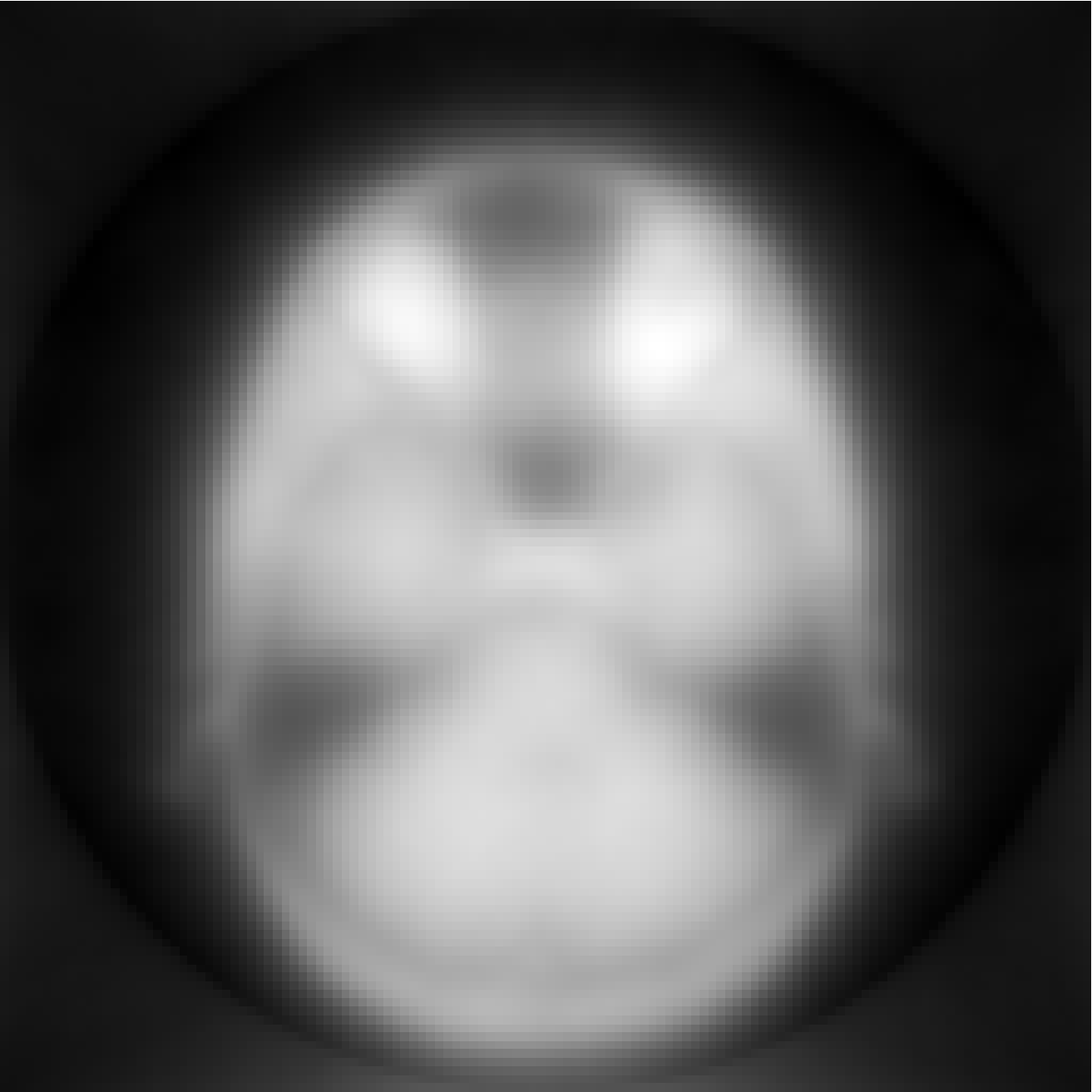}}   &
    {\label{Fig:mri ABGKB GCV RNS}\includegraphics[width=.18\textwidth]{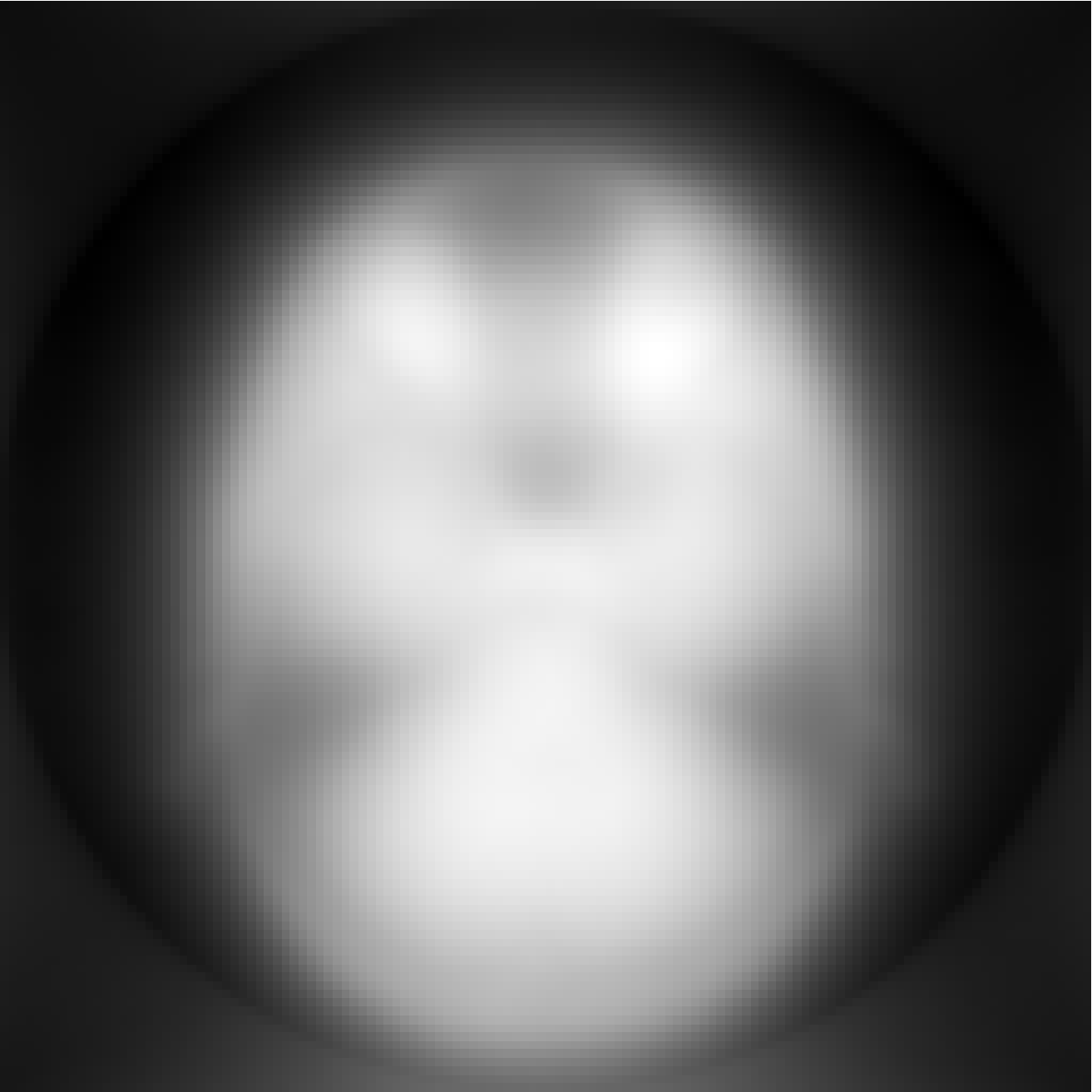}} &
    {\label{Fig:mri ABGKB Lcurve DP}\includegraphics[width=.18\textwidth]{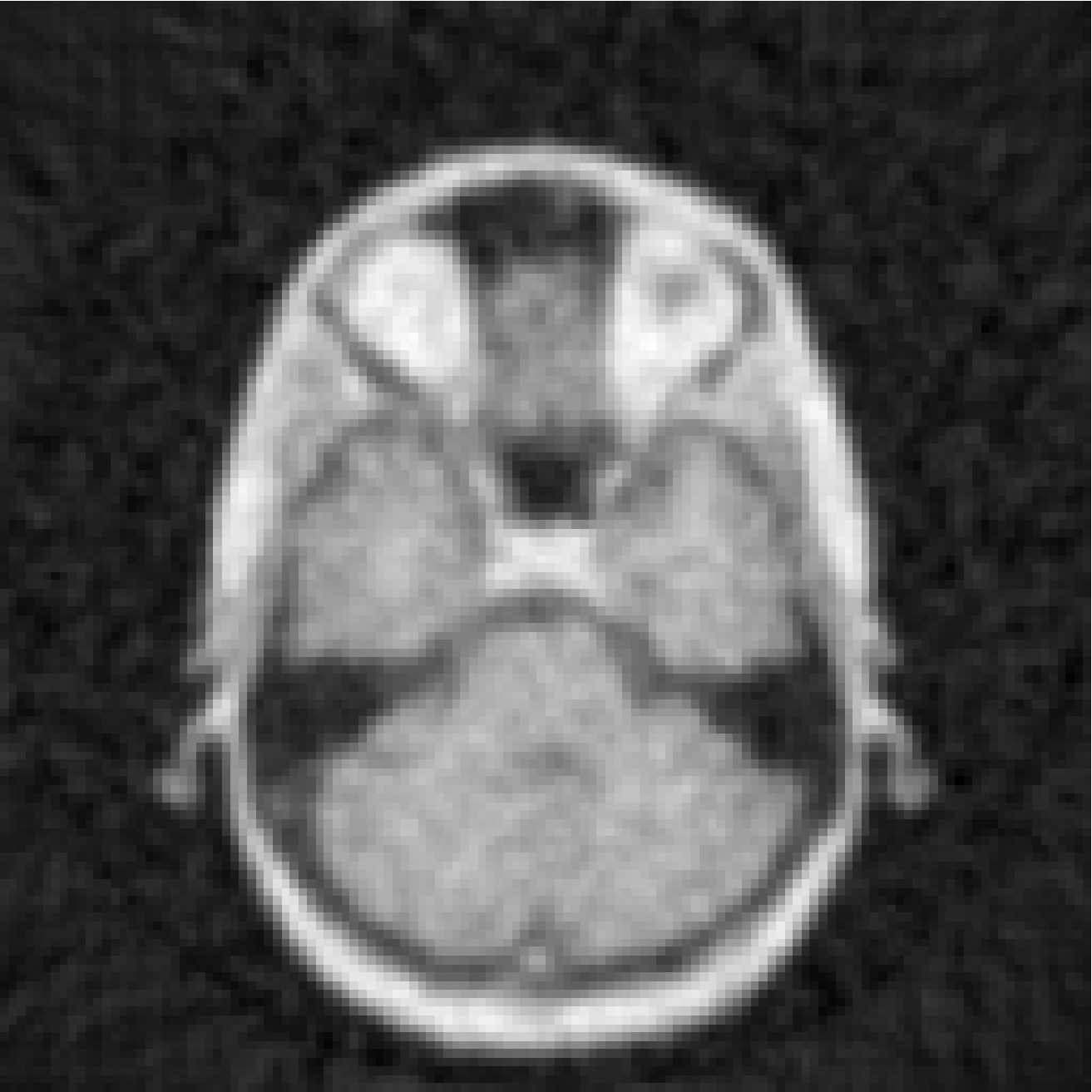}}   &
     {\label{Fig:mri ABGKB Lcurve RNS}\includegraphics[width=.18\textwidth]{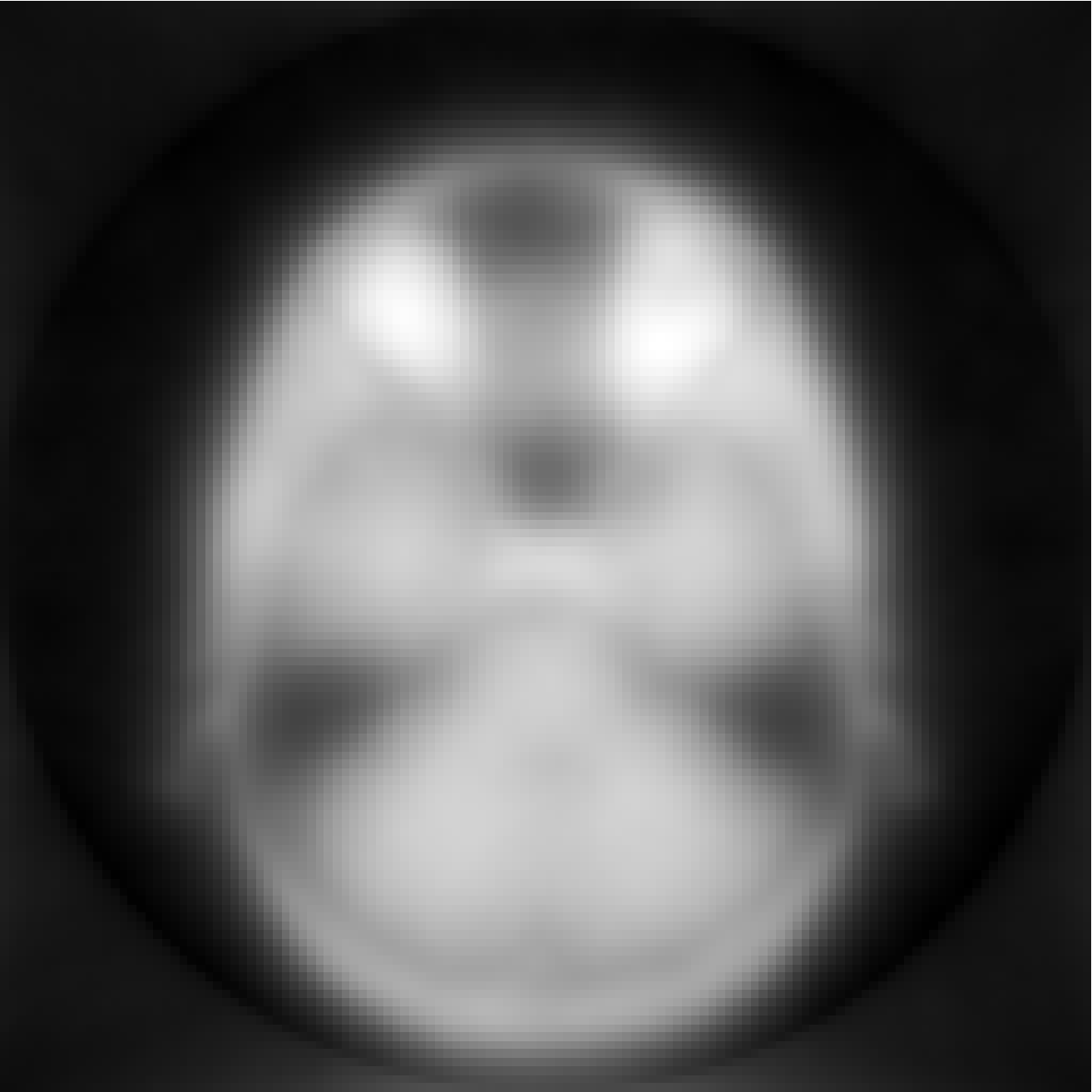}}   \\ 
      {\label{Fig:mri BAGKB GCV DP}\includegraphics[width=.18\textwidth]{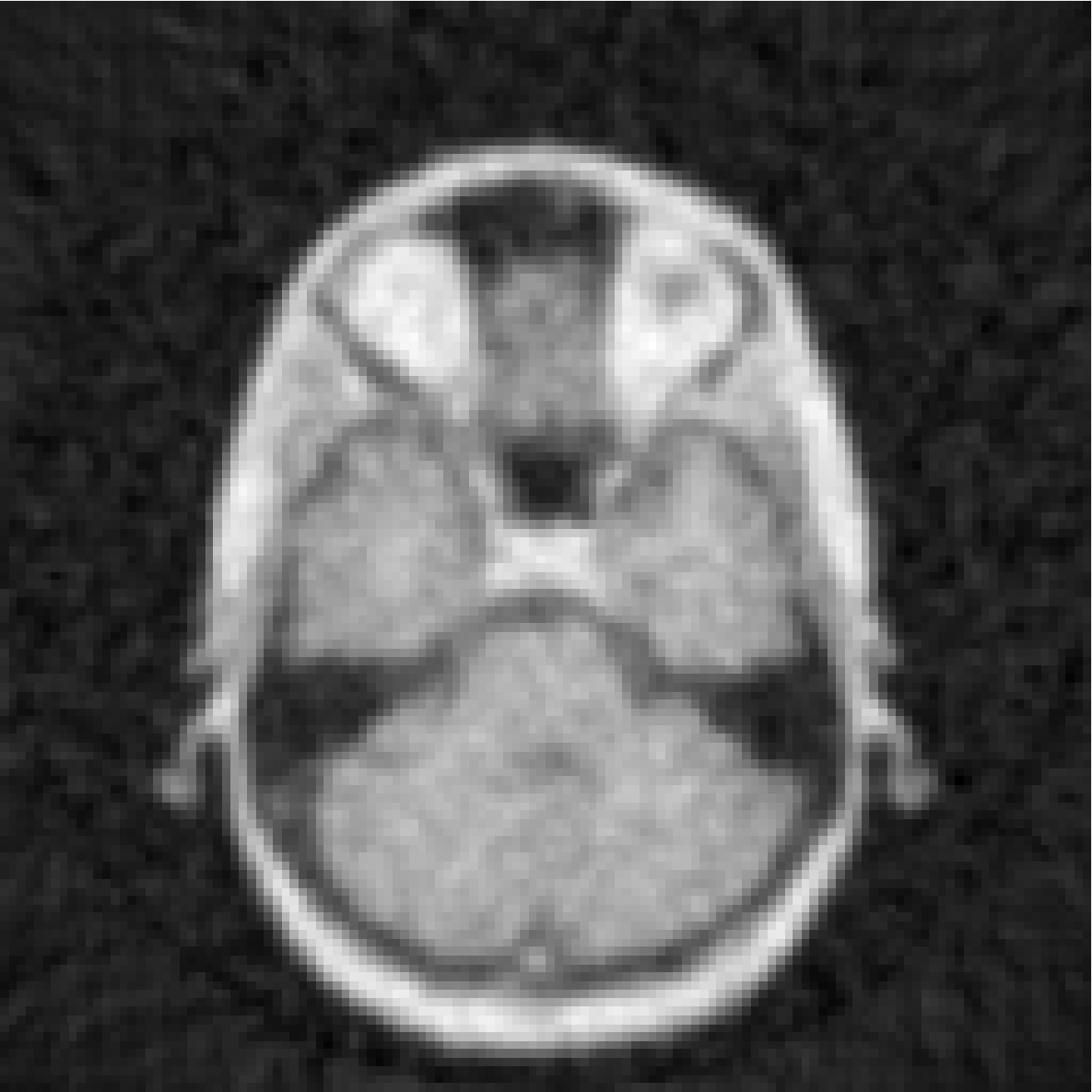}} &
      {\label{Fig:mri BAGKB GCV RNS}\includegraphics[width=.18\textwidth]{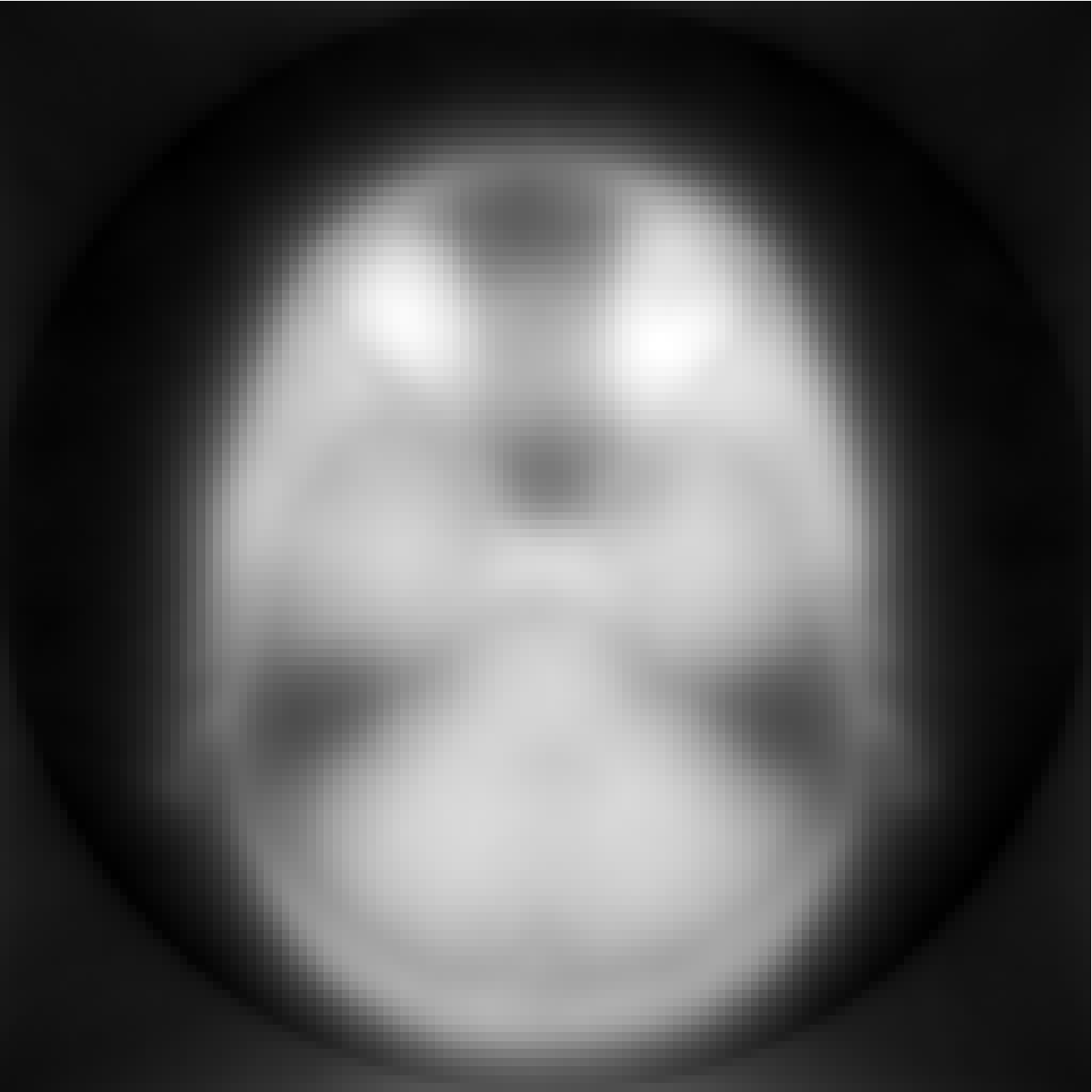}} &
     {\label{Fig:mri BAGKB Lcurve DP}\includegraphics[width=.18\textwidth]{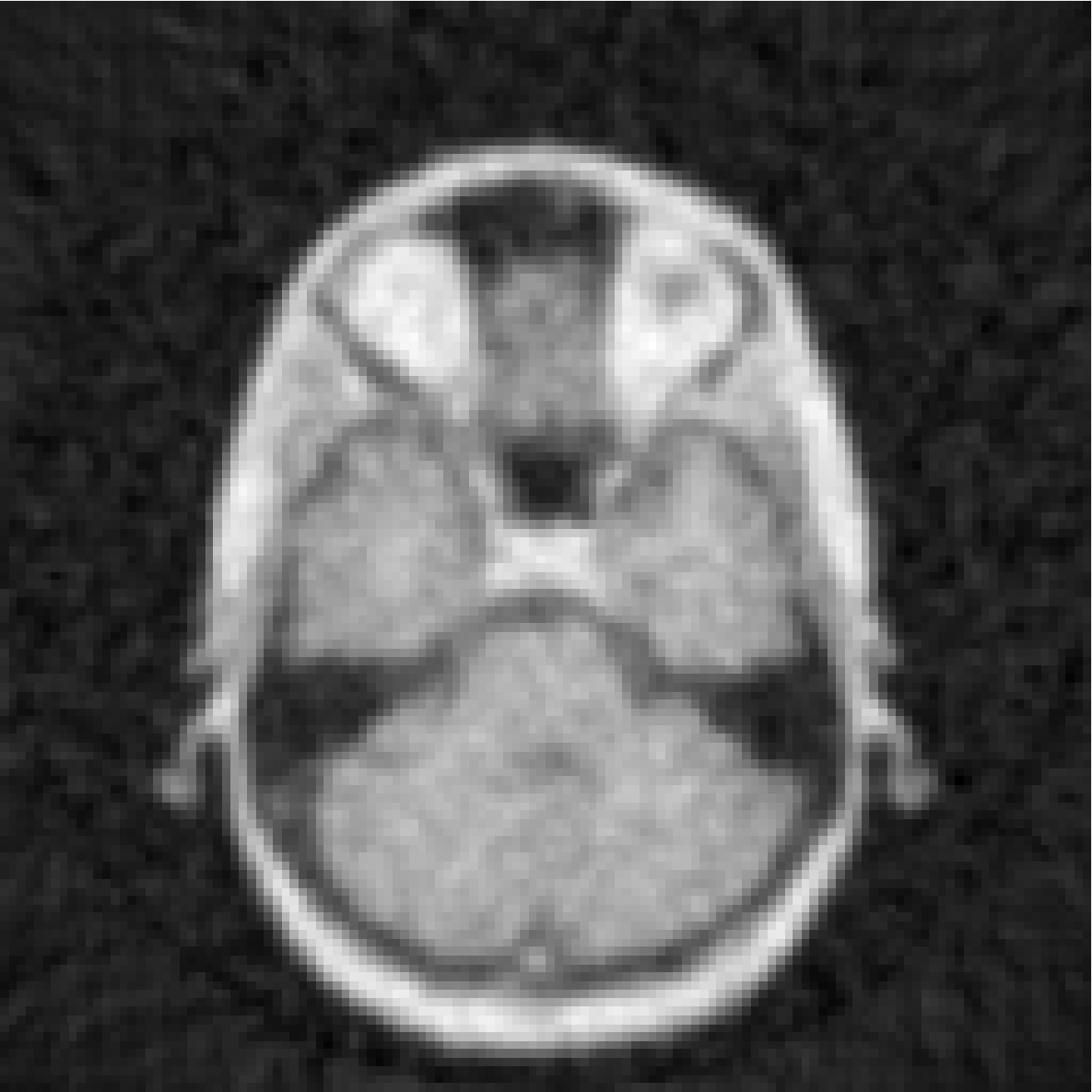}} &
      {\label{Fig:mri BAGKB Lcurve RNS}\includegraphics[width=.18\textwidth]{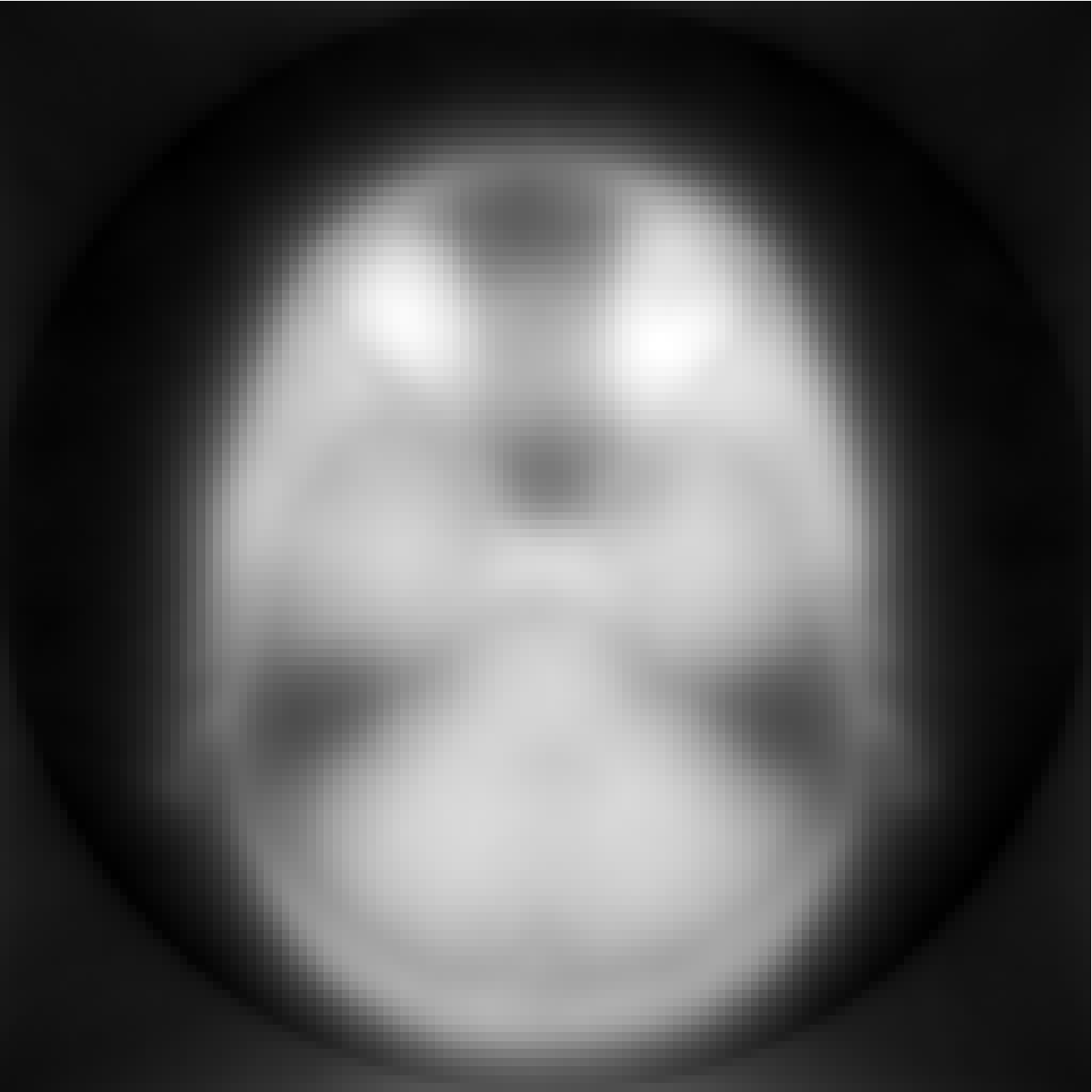}} \\
 {\label{Fig:mri ABGMRES GCV DP}\includegraphics[width=.18\textwidth]{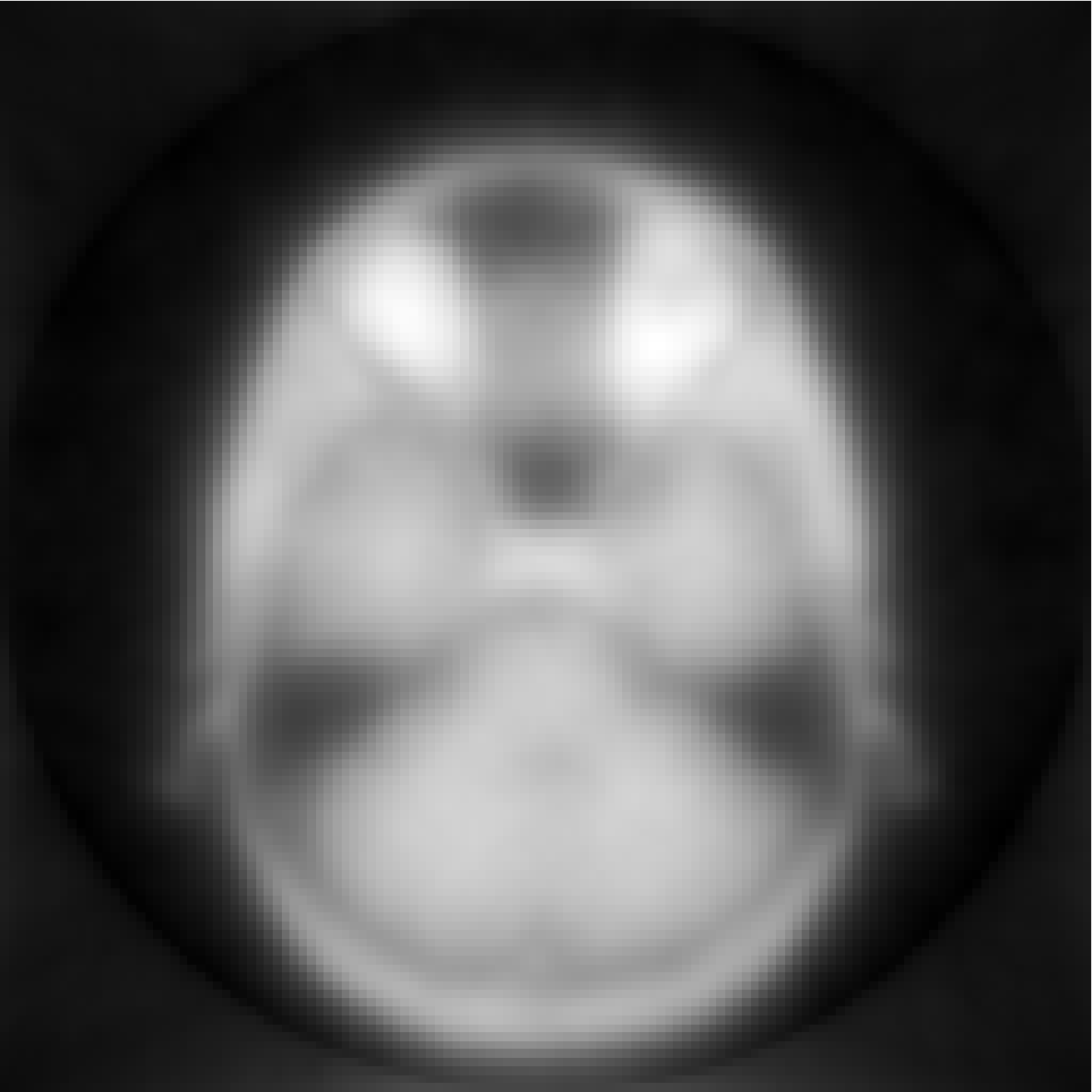}} & 
    {\label{Fig:mri ABGMRES GCV RNS}\includegraphics[width=.18\textwidth]{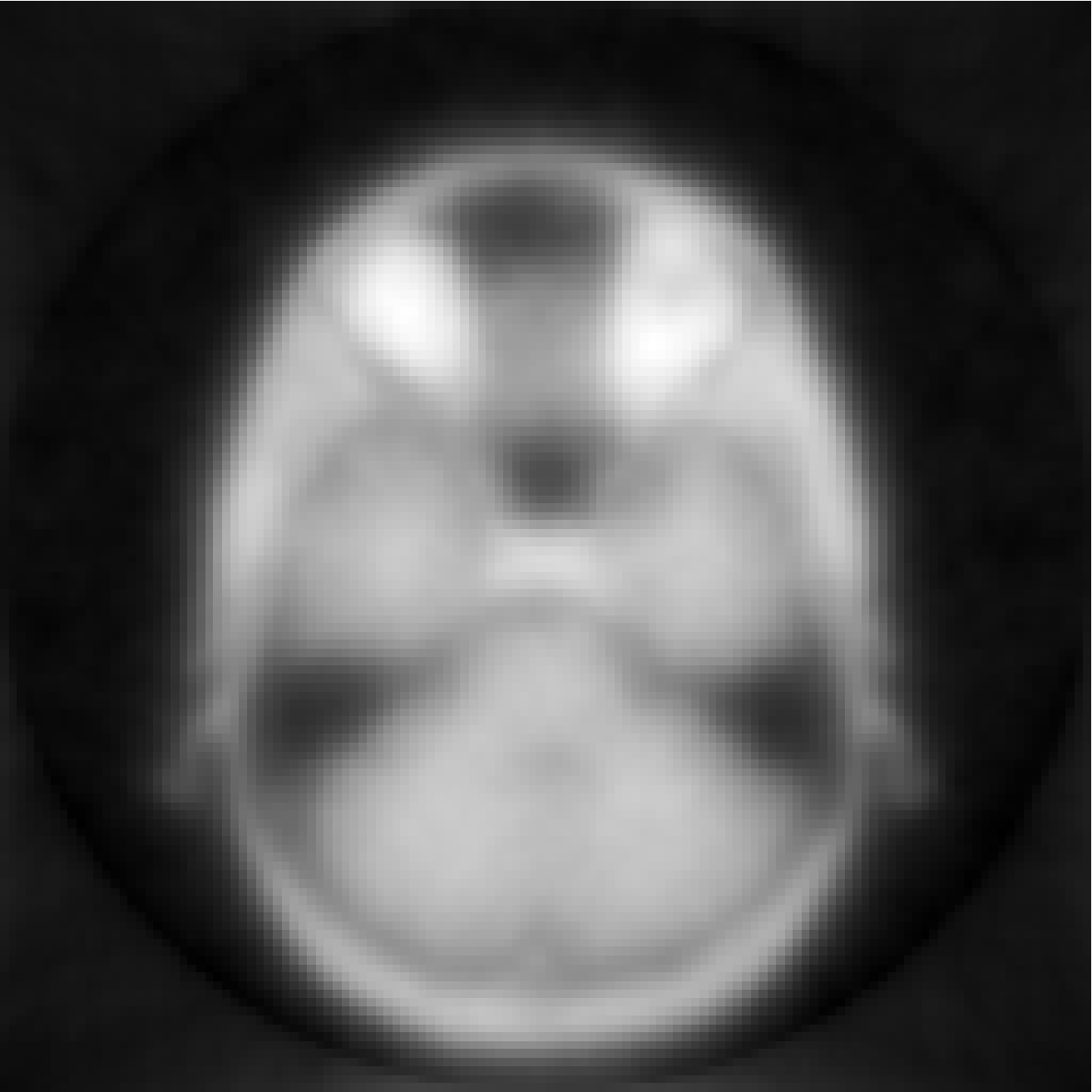}}  &
   {\label{Fig:mri ABGMRES Lcurve DP}\includegraphics[width=.18\textwidth]{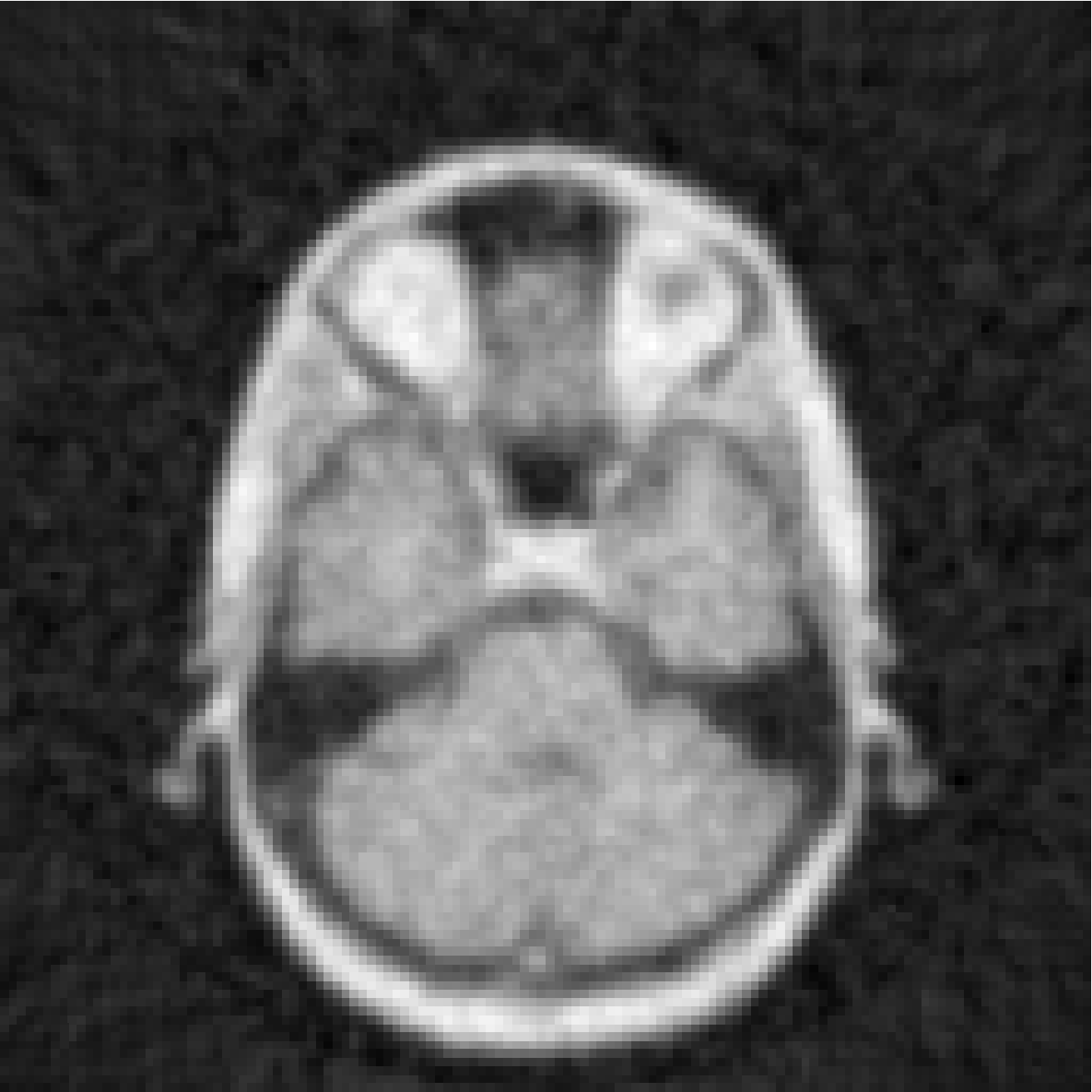}}  &
   {\label{Fig:mri ABGMRES Lcurve RNS}\includegraphics[width=.18\textwidth]{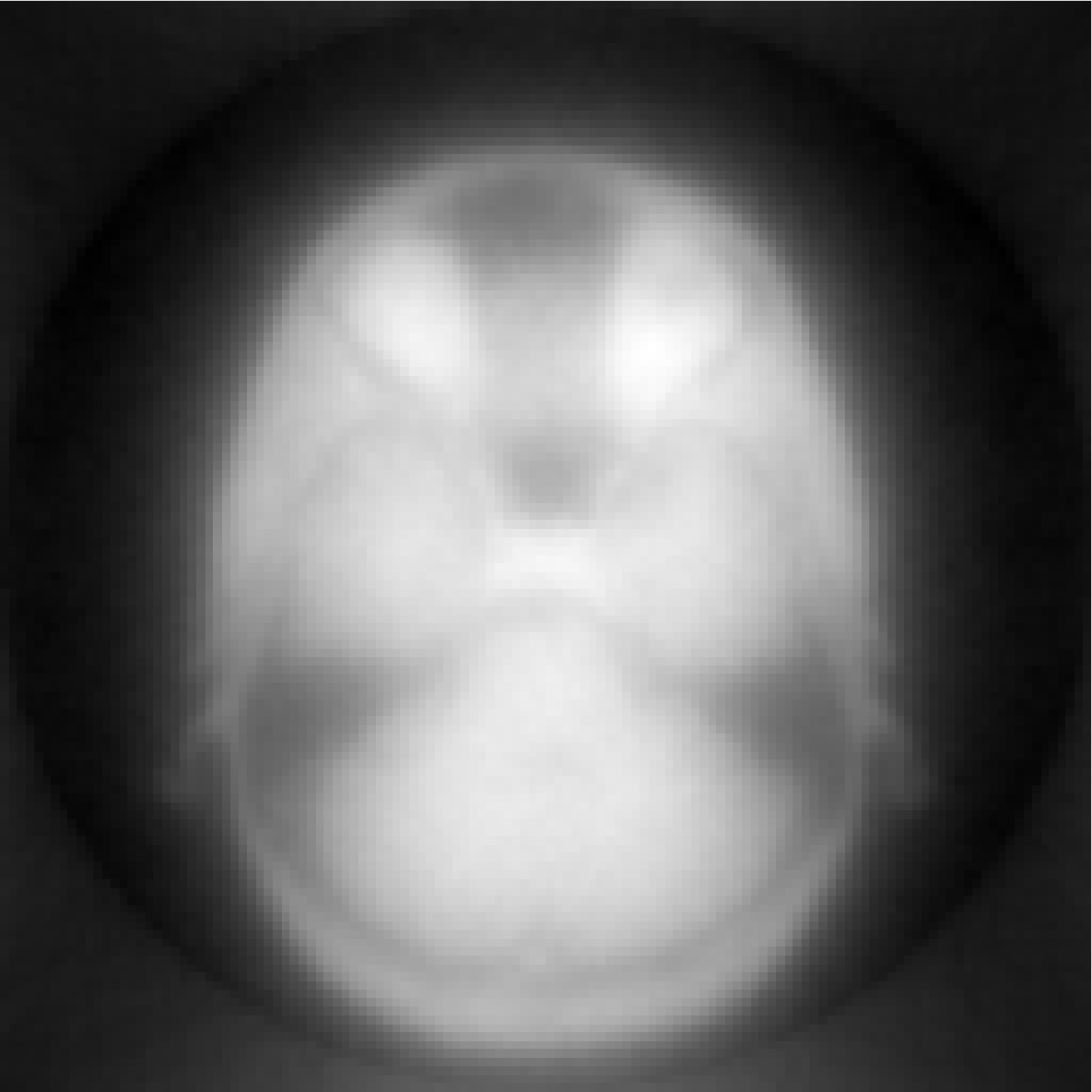}}  \\  
  {\label{Fig:mri BAGMRES GCV DP}\includegraphics[width=.18\textwidth]{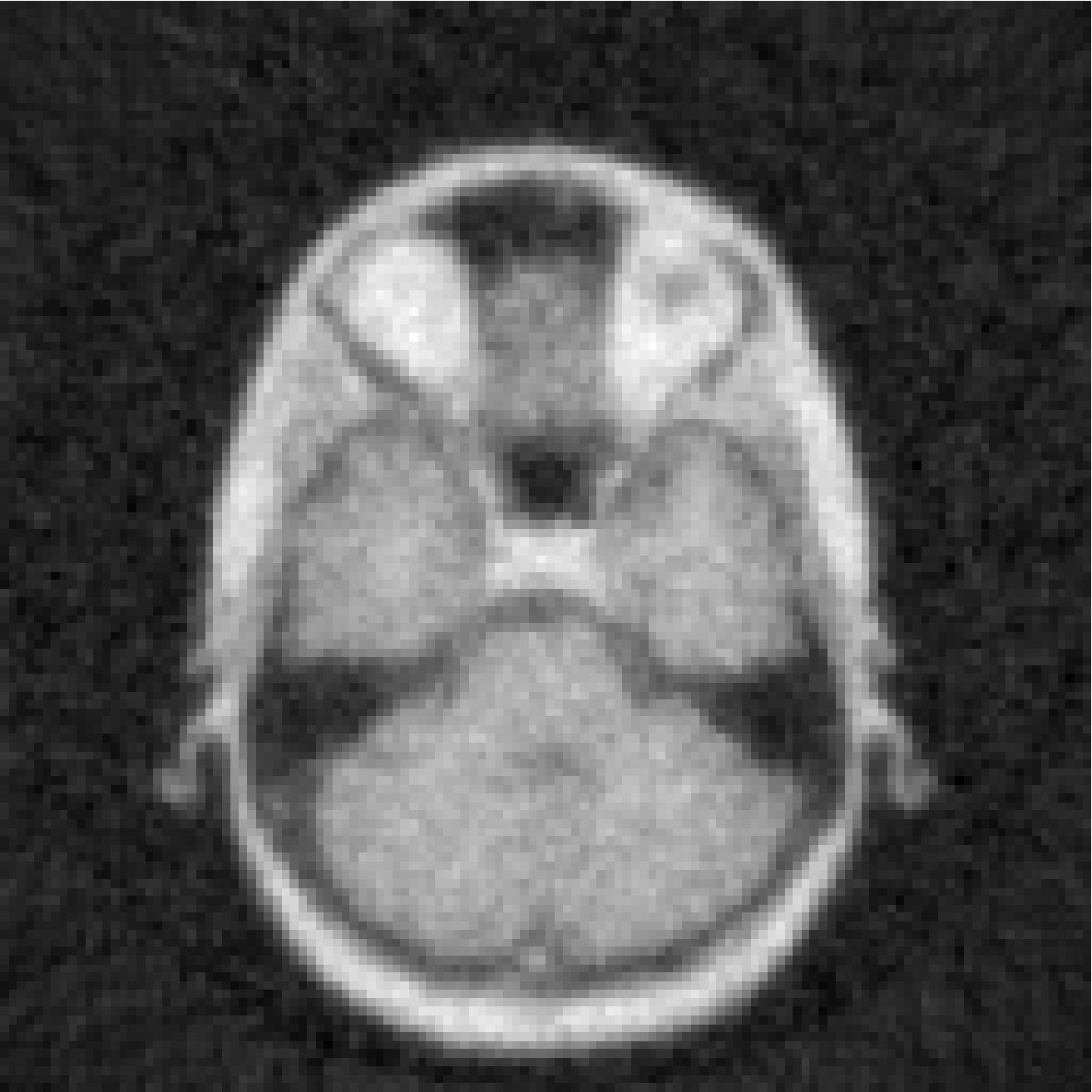}}  &{\label{Fig:mri BAGMRES GCV RNS}\includegraphics[width=.18\textwidth]{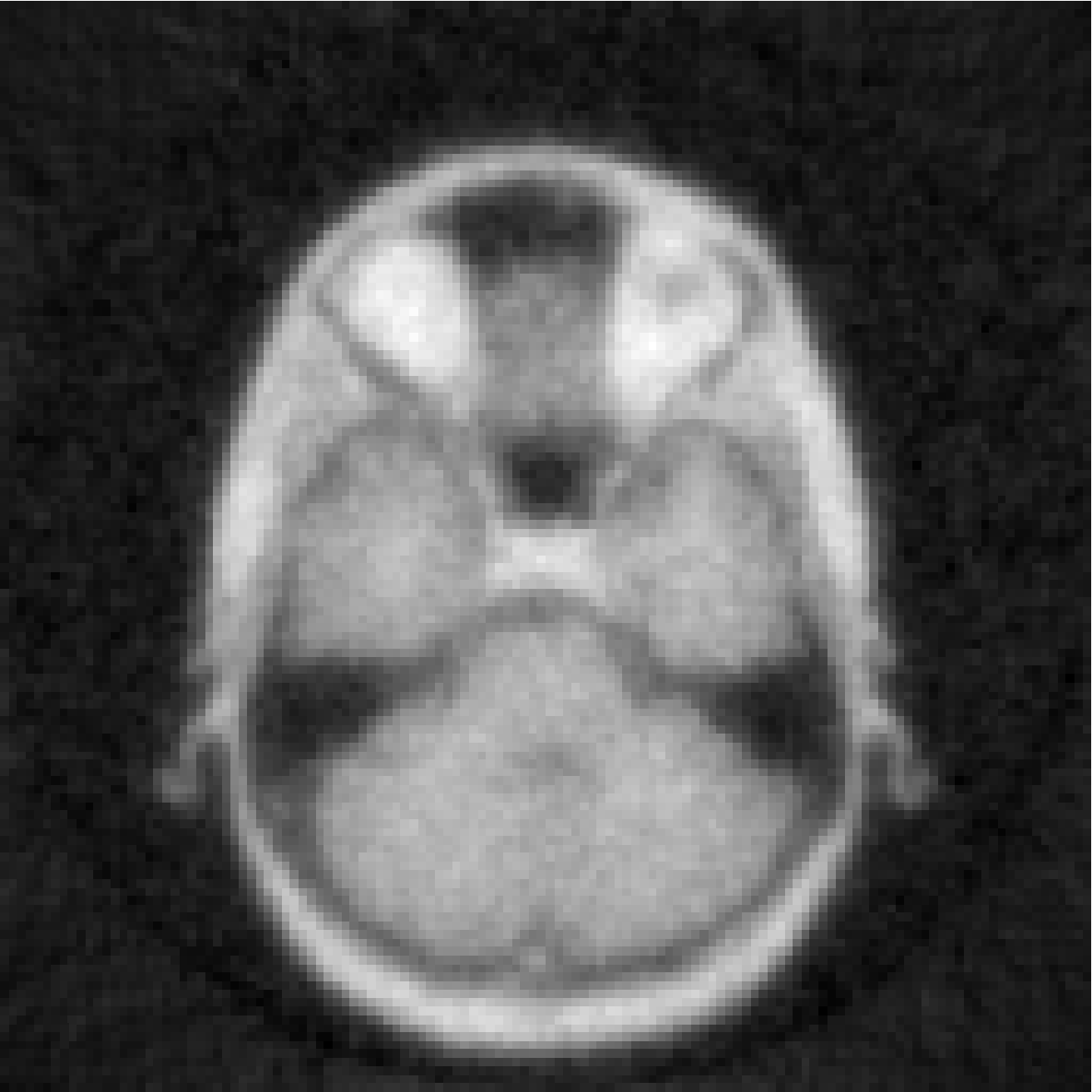}}  &{\label{Fig:mri BAGMRES Lcurve DP}\includegraphics[width=.18\textwidth]{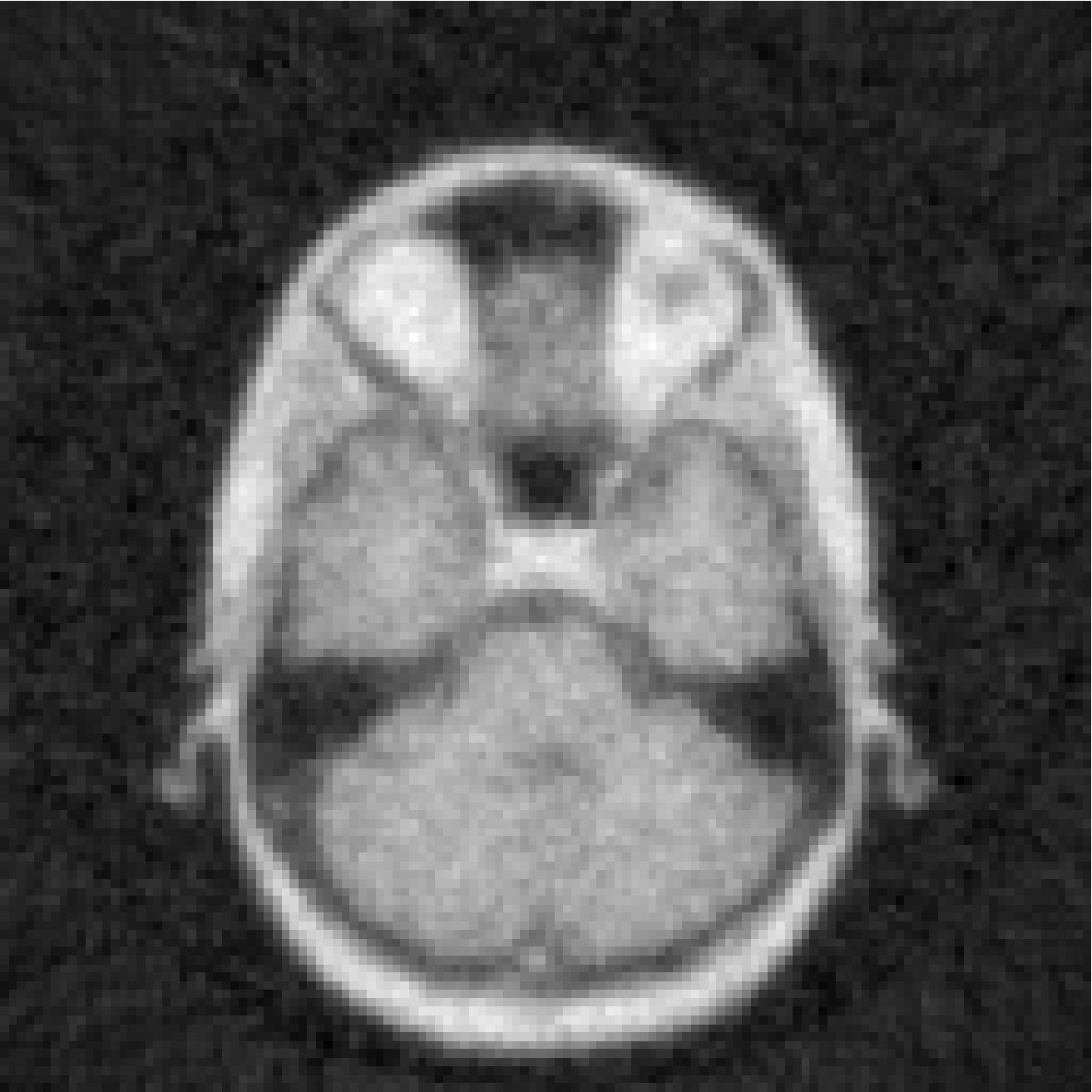}}  &{\label{Fig:mri BAGMRES Lcurve RNS}\includegraphics[width=.18\textwidth]{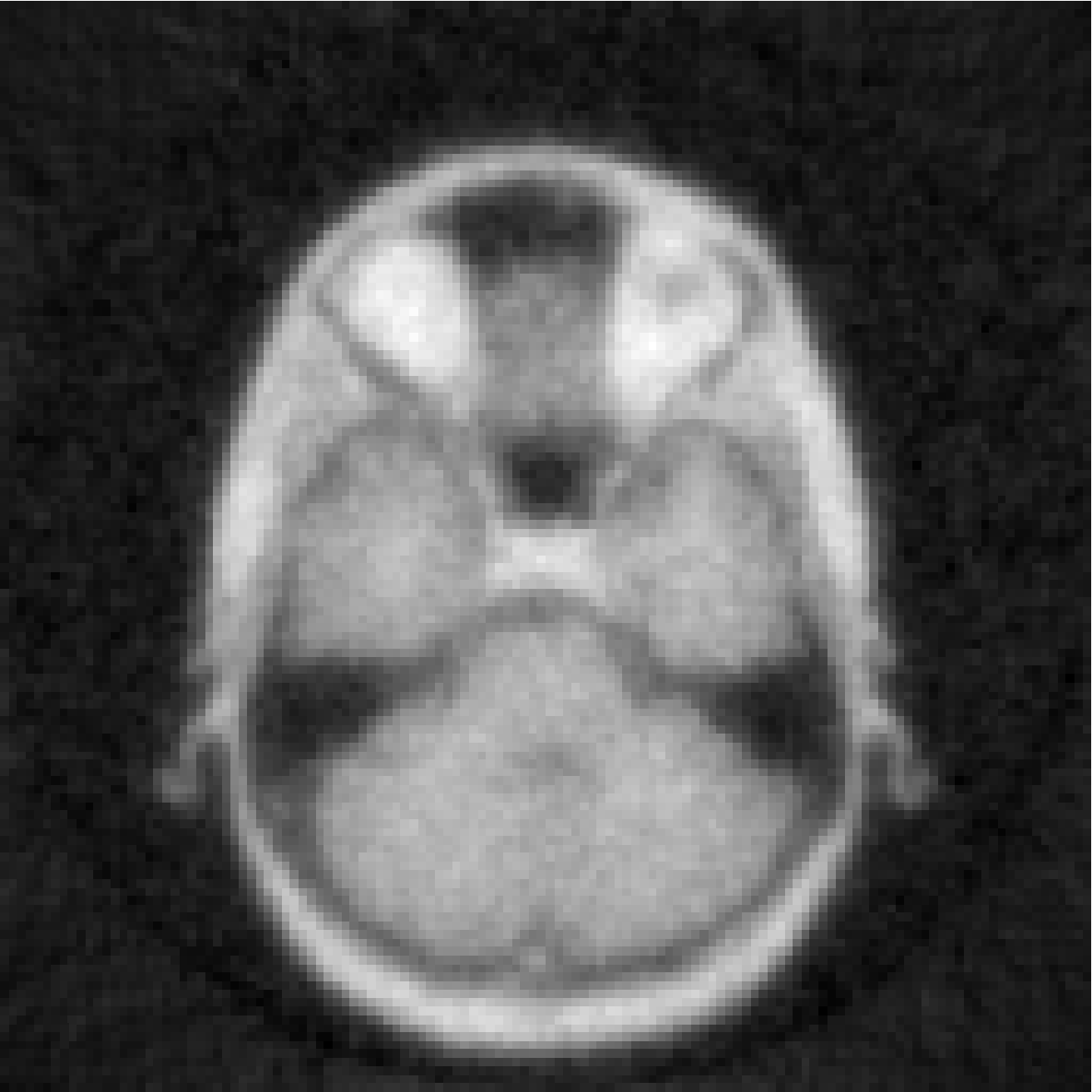}} 
  \end{tabular}
  \caption{Reconstructed images corresponding to \Cref{Fig:blurred mri}, obtained using hybrid AB-GKB, hybrid BA-GKB, hybrid AB-GMRES, and hybrid BA-GMRES with methods in columns, shown in the first, second, third, and fourth rows, respectively.}\label{Fig:Recons mri} 
  \end{figure}

\begin{example}\label{ex:problem1}[Overdetermined Example]
We consider a phantom image of size $128 \times 128$, shown in \Cref{Fig:True mri}. The sinogram data are generated using parallel geometry through $180$ view angles. We use the ASTRA Toolbox \citet{van2015astra} to generate the forward projector $A \in \mathbb{R}^{23040\times 16384}$  and backprojector $B \in \mathbb{R}^{16384\times 23040}$.
Then we use $6 \%$ white Gaussian noise, corresponding to SNR of approximately $24$, to obtain the contaminated sinogram data $\bfb \in \mathbb{R}^{23040}$, shown in \Cref{Fig:blurred mri}. The convergence results are illustrated in \Cref{Fig:mri RRE GCV,Fig:mri RRE Lcurve} for the GCV and L-curve, respectively, and quantitative results are reported in \Cref{Tab:Small Problem}. The corresponding reconstructed images are presented in \Cref{Fig:Recons mri}. 
\end{example}

\begin{table}
 \footnotesize 
 \caption{\label{Tab:Small Problem} The corresponding \text{RRE} \eqref{eq:RE} computed at $k_{DP}$ \eqref{eq:DP} and $k_{RNS}$ \eqref{eq:RNS} in parentheses, by all methods with $\lambda$ obtained with the GCV and L-curve, applied to the problem with $\text{SNR}\approx 24$ in \Cref{ex:problem1}. The lowest RRE results are shown in boldface. The timings are for all algorithms out to $150$ iterations and averaged over five runs.}
 \begin{tabular}[t]{lccccccccc } 
 \toprule 
Method&Regularization&Iteration&$\text{RRE}(\bfx)$&Time(s)  \\ \midrule
\multirow{2}{*}{hybrid AB-GKB} & GCV& $150(5)$ &$0.42(0.49)$& \multirow{2}{*}{$1.82$}  \\ 
& L-curve & $26(7)$& $\boldsymbol{0.20}(0.40)$ &  \\ 
\hline
\multirow{2}{*}{hybrid BA-GKB} & GCV& $31(7)$ &$\boldsymbol{0.20}(0.41)$& \multirow{2}{*}{$1.79$}  \\ 
& L-curve & $31(7)$& $\boldsymbol{0.20}(0.41)$ &  \\ 
\hline
\multirow{2}{*}{hybrid AB-GMRES} & GCV& $150(7)$ &$0.38(0.35)$& \multirow{2}{*}{$1.34$}  \\ 
& L-curve& $12(4)$& $\boldsymbol{0.20}(0.45)$ &  \\ 
\hline
\multirow{2}{*}{hybrid BA-GMRES} & GCV& $10(8)$ &$0.21(0.24)$& \multirow{2}{*}{$1.26$}  \\ 
& L-curve & $10(8)$& $0.21(0.24)$ & \\ 
\bottomrule
 \end{tabular}    
 \end{table}

From \Cref{Fig:mri RRE GCV}, it is immediately apparent that the hybrid BA-GMRES obtained using GCV is significantly affected by semiconvergence compared to all other methods. Furthermore, we can see that hybrid BA-GKB works well with GCV and yields excellent performance. The plots in \Cref{Fig:mri RRE Lcurve} illustrate that the L-curve offers more stable solutions that are comparable in RRE to those derived from GCV. However, the hybrid BA-GMRES with the L-curve still exhibits semiconvergence, but only for a higher number of iterations than when using the GCV.  The results presented in column $4$ of \Cref{Tab:Small Problem} support the effectiveness of DP against the RNS in terminating algorithms. Furthermore, we observe that while GCV and L-curve significantly affect the performance of the hybrid AB-GKB and AB-GMRES, this is not the case for their BA counterparts. We should emphasize that although $\epsilon$ in the RNS is selected to minimize the GMRES methods, the GKB methods with DP can yield comparable RRE.

The execution times listed in column $5$ of \Cref{Tab:Small Problem} show that the GMRES methods are faster than the GKB methods. Furthermore, the BA algorithms are faster than their AB counterparts. This validates the computational analysis presented in \Cref{sec:comp costs}, \crefrange{cost:ABGMRES}{cost:BAGKB}.

\begin{figure}
  \centering
   \subfloat[Slice 5]{\label{Fig:True slice 5}\includegraphics[width=.25\textwidth]{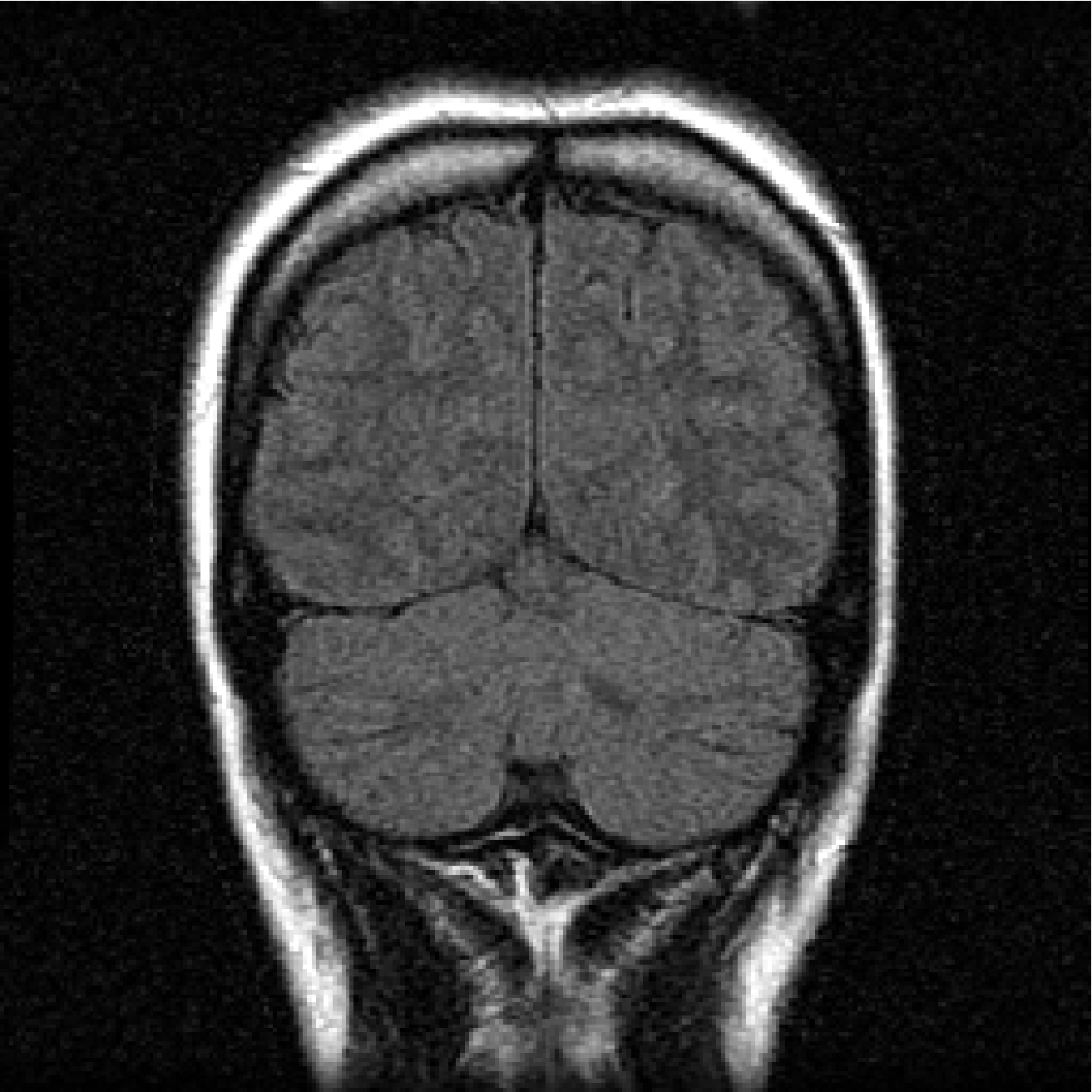}} \
       \subfloat[Slice 11]{\label{Fig:True slice 11}\includegraphics[width=.25\textwidth]{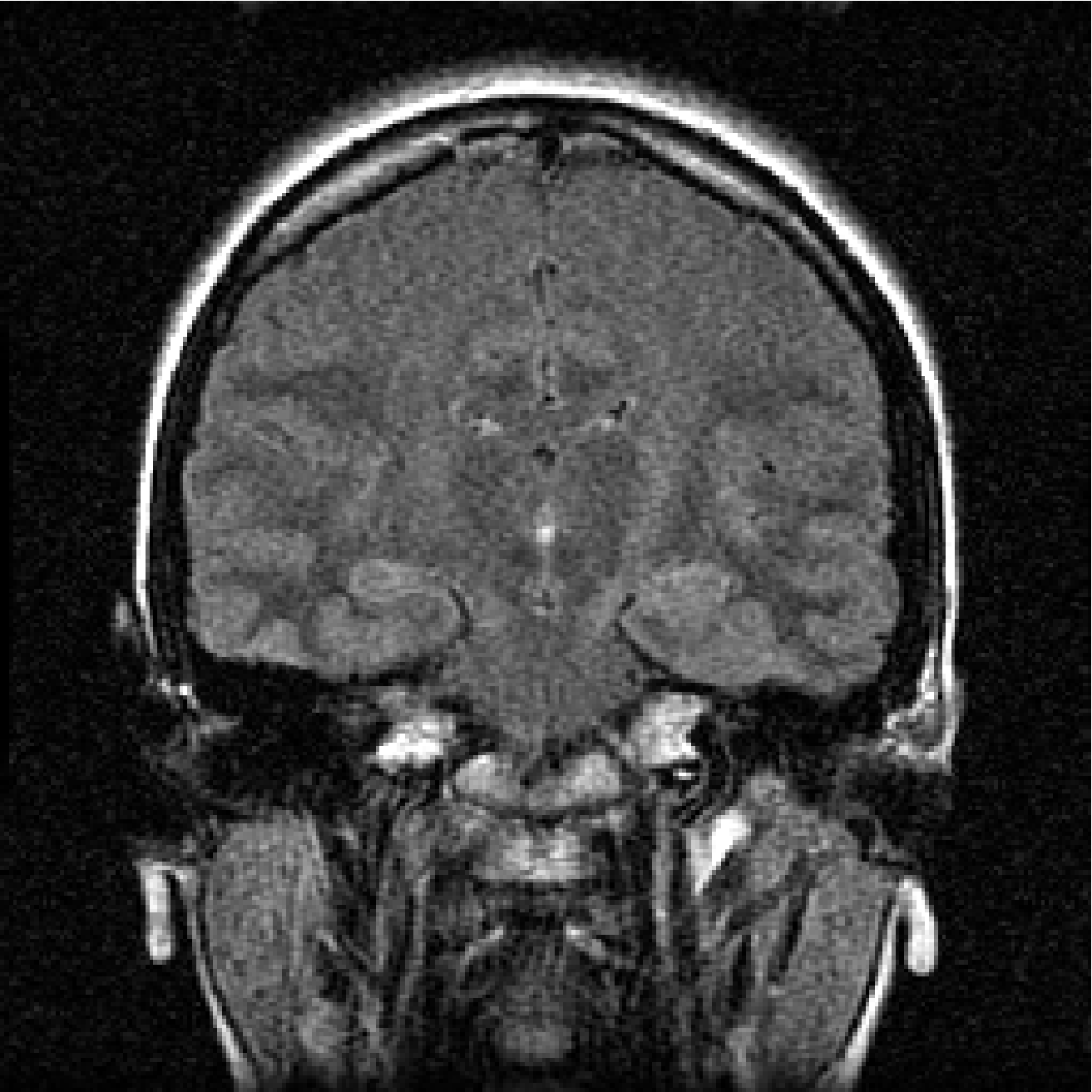}} \
    \subfloat[Slice 18]{\label{Fig:True slice 18}\includegraphics[width=.25\textwidth]{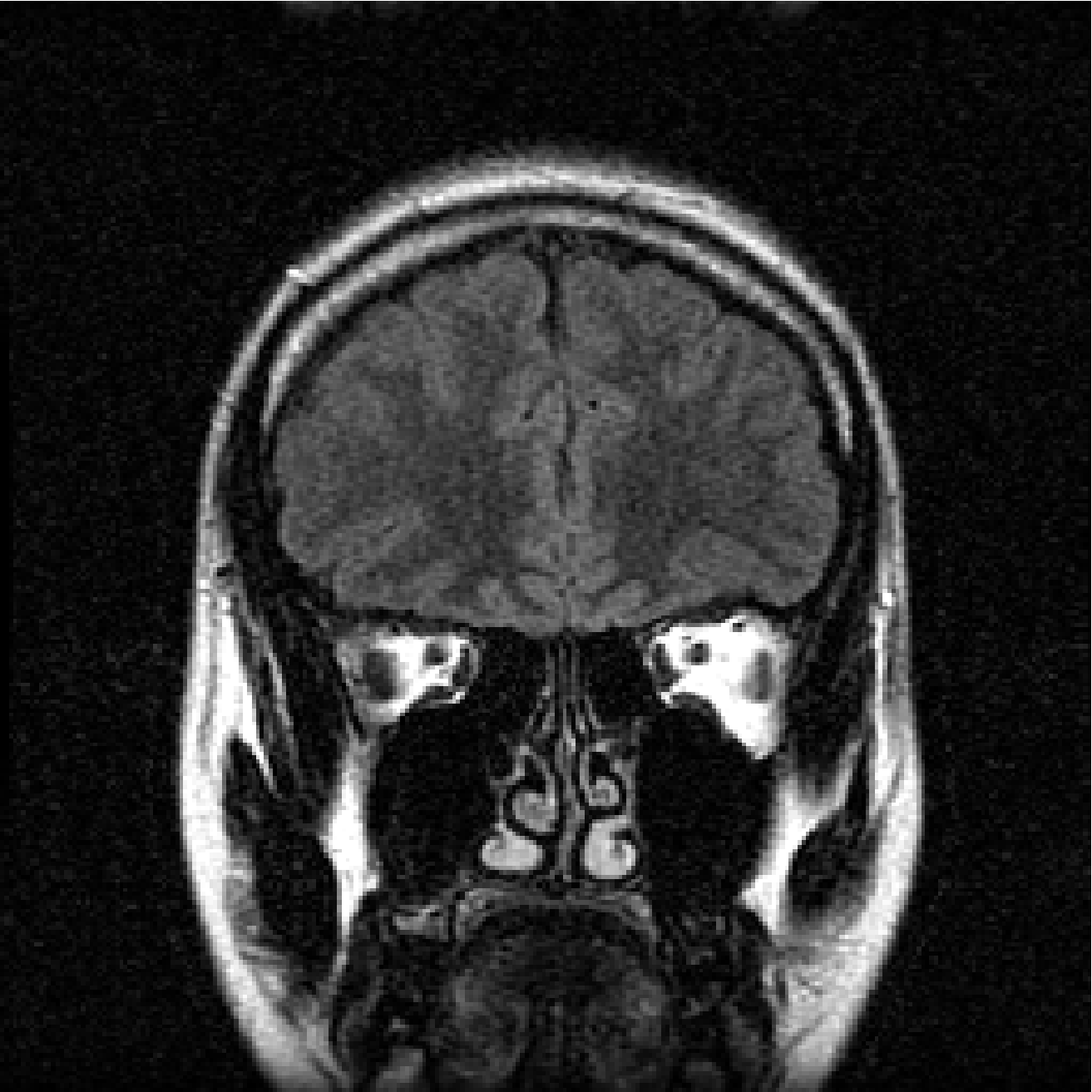}} \\
   \subfloat[Slice 5]{\label{Fig:blurred slice 5}\includegraphics[width=.25\textwidth]{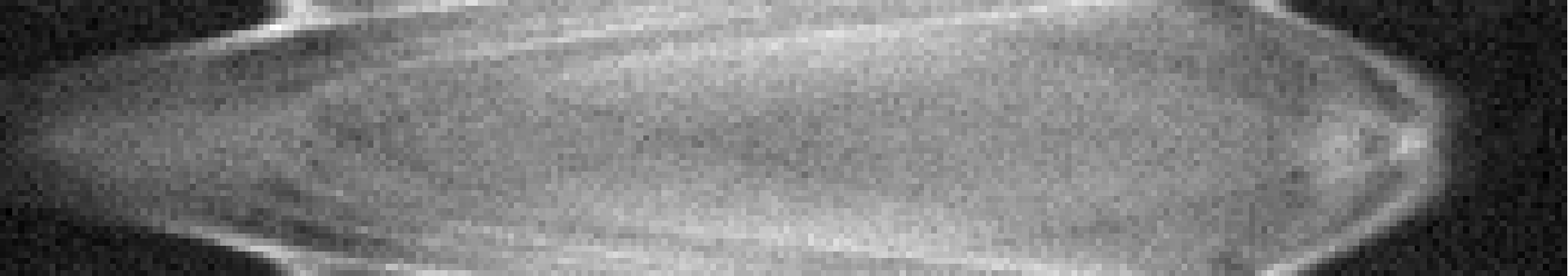}}\
   \subfloat[Slice 11]{\label{Fig:blurred slice 11}\includegraphics[width=.25\textwidth]{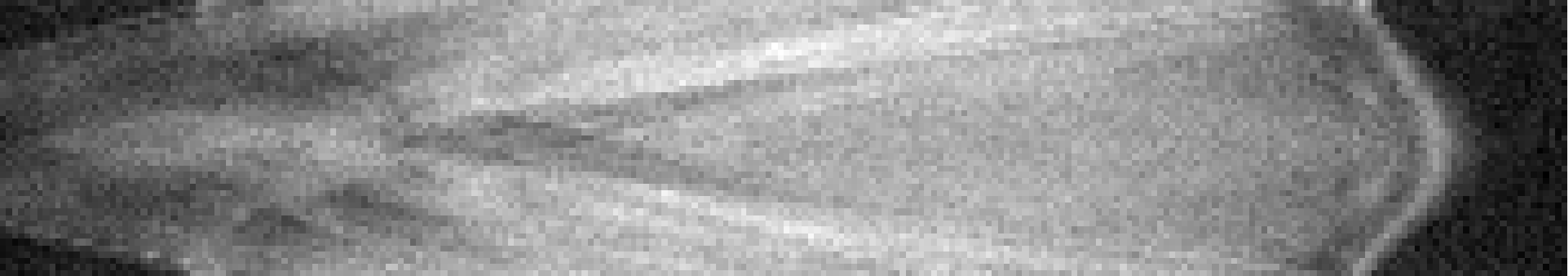}}\
   \subfloat[Slice 18]{\label{Fig:blurred slice 18}\includegraphics[width=.25\textwidth]{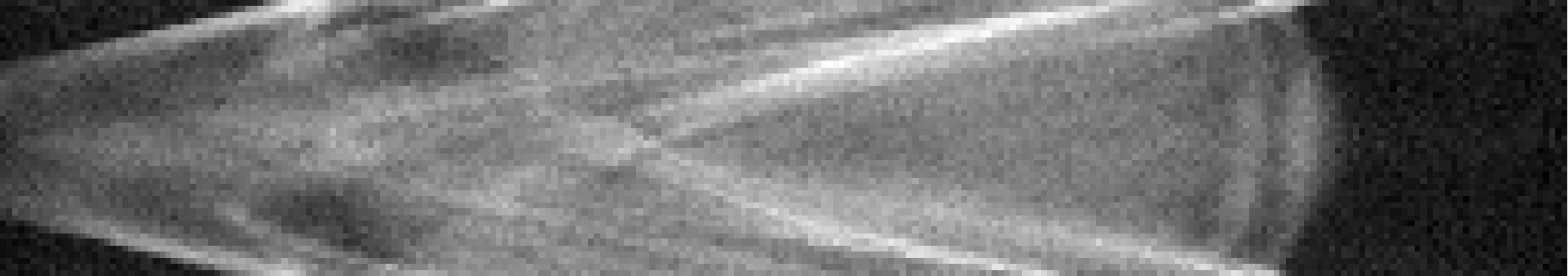}}
 \caption{The true images of $256 \times 256$ pixels in the first row and the corresponding noisy sinograms with $8\%$ Gaussian noise in the second row, for the problem in \Cref{ex:problem2}. \label{Fig:slices True and noisy}}  
  \end{figure}

\begin{figure}
  \centering
   \subfloat[GCV]{\label{Fig:slices RRE GCV}\includegraphics[width=.48\textwidth]{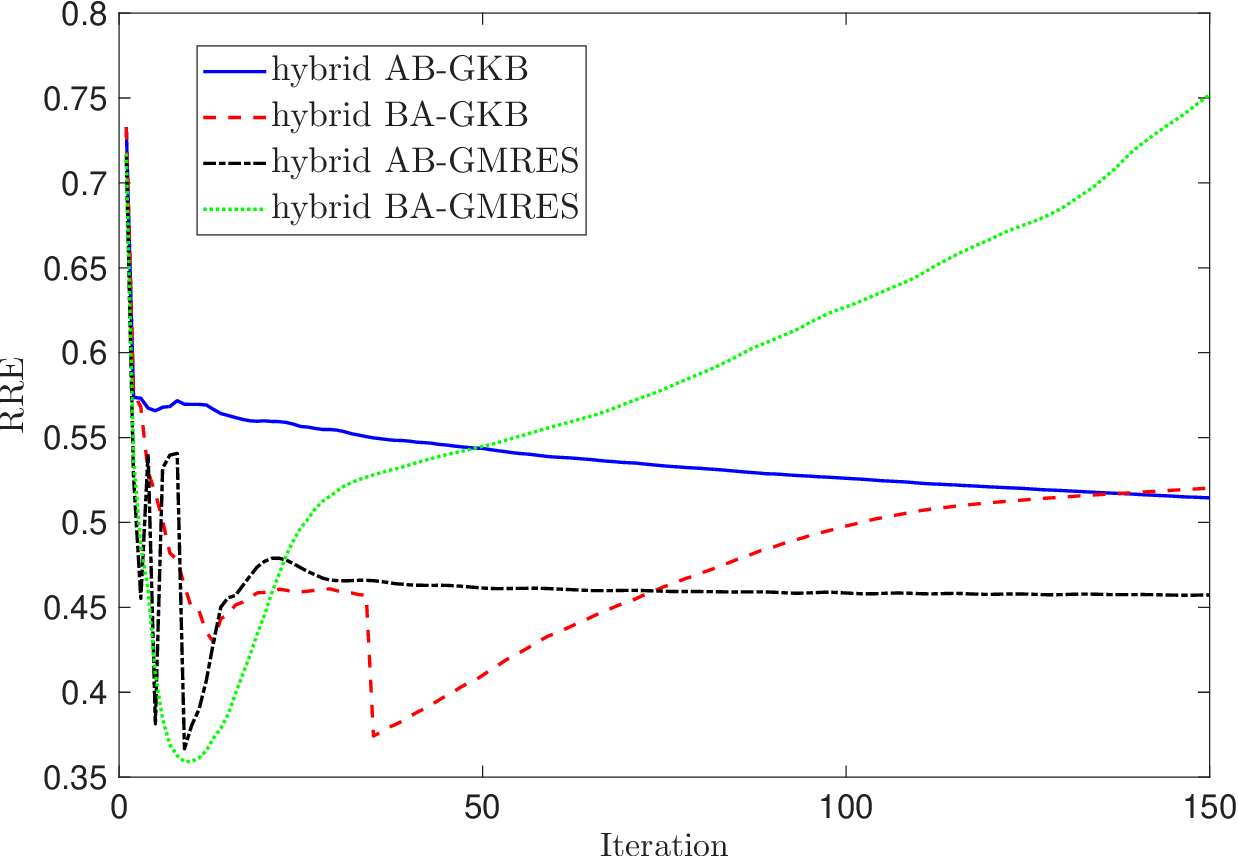}}   \ 
  \subfloat[L-curve]{\label{Fig:slices RRE Lcurve}\includegraphics[width=.48\textwidth]{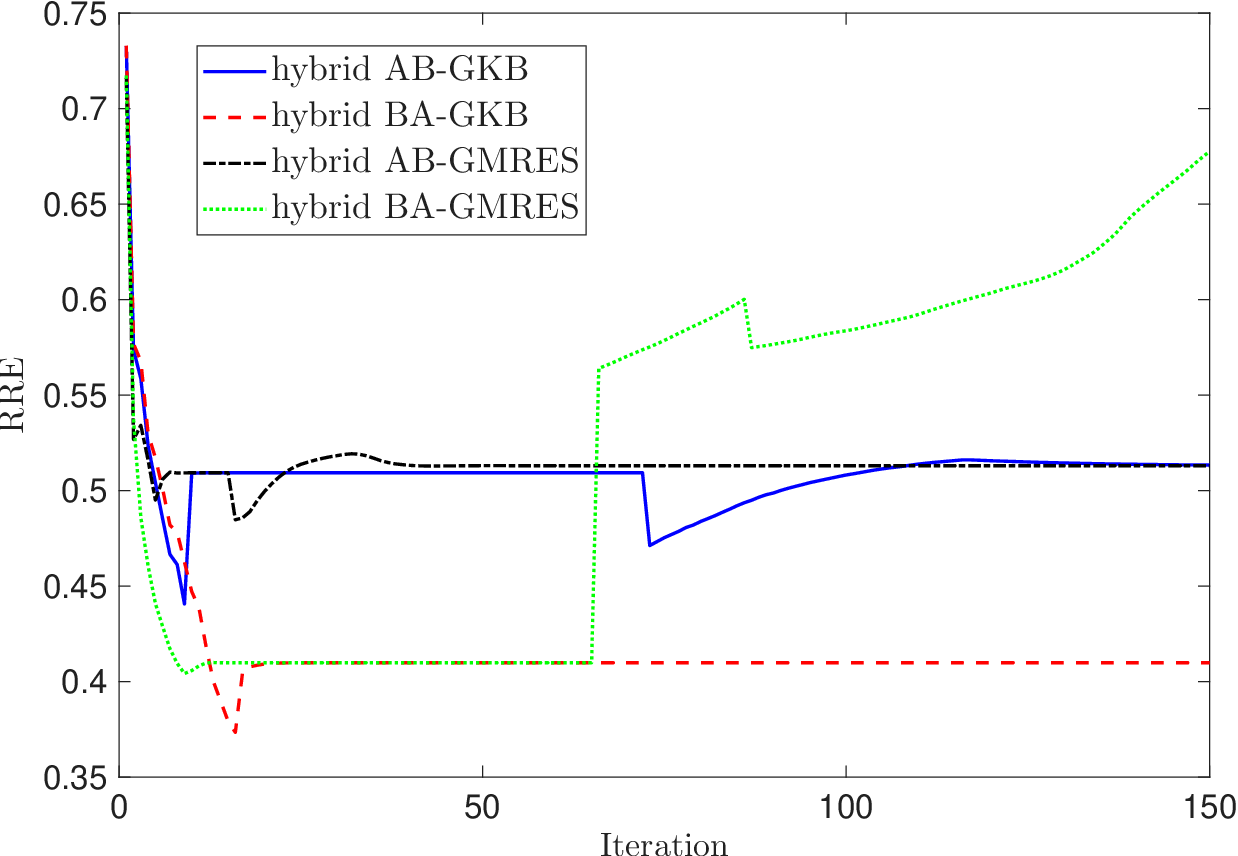}} 
  \caption{RRE for slice $11$ with different methods for unmatched projectors  in \Cref{ex:problem2}.}\label{Fig:slices}  
  \end{figure}

\begin{figure}
  \centering
   \begin{tabular}{ccc} 
 {\label{Fig:ABGKBGCVDP 5}\includegraphics[width=.20\textwidth]{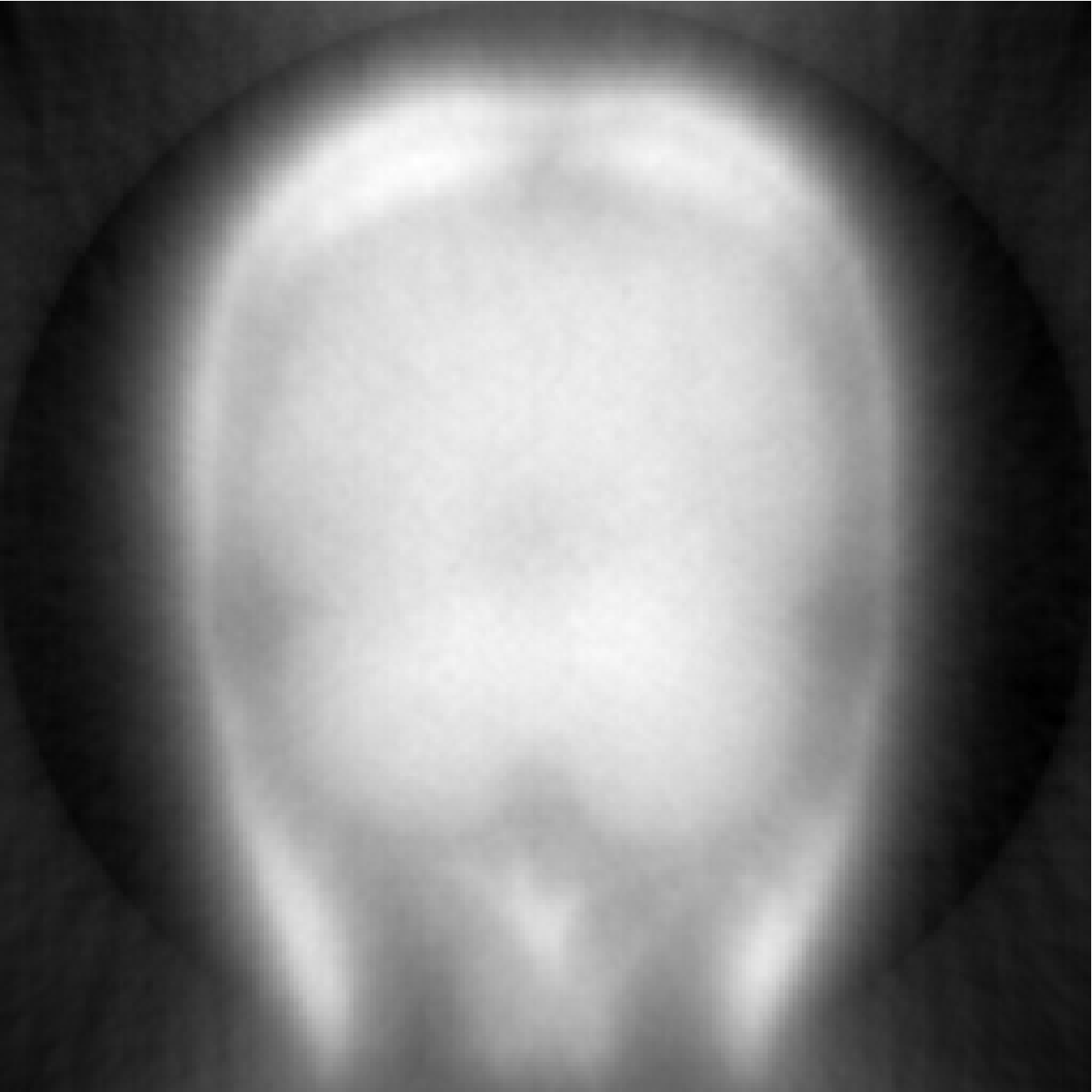}} &
 {\label{Fig:ABGKBGCVDP 11}\includegraphics[width=.20\textwidth]{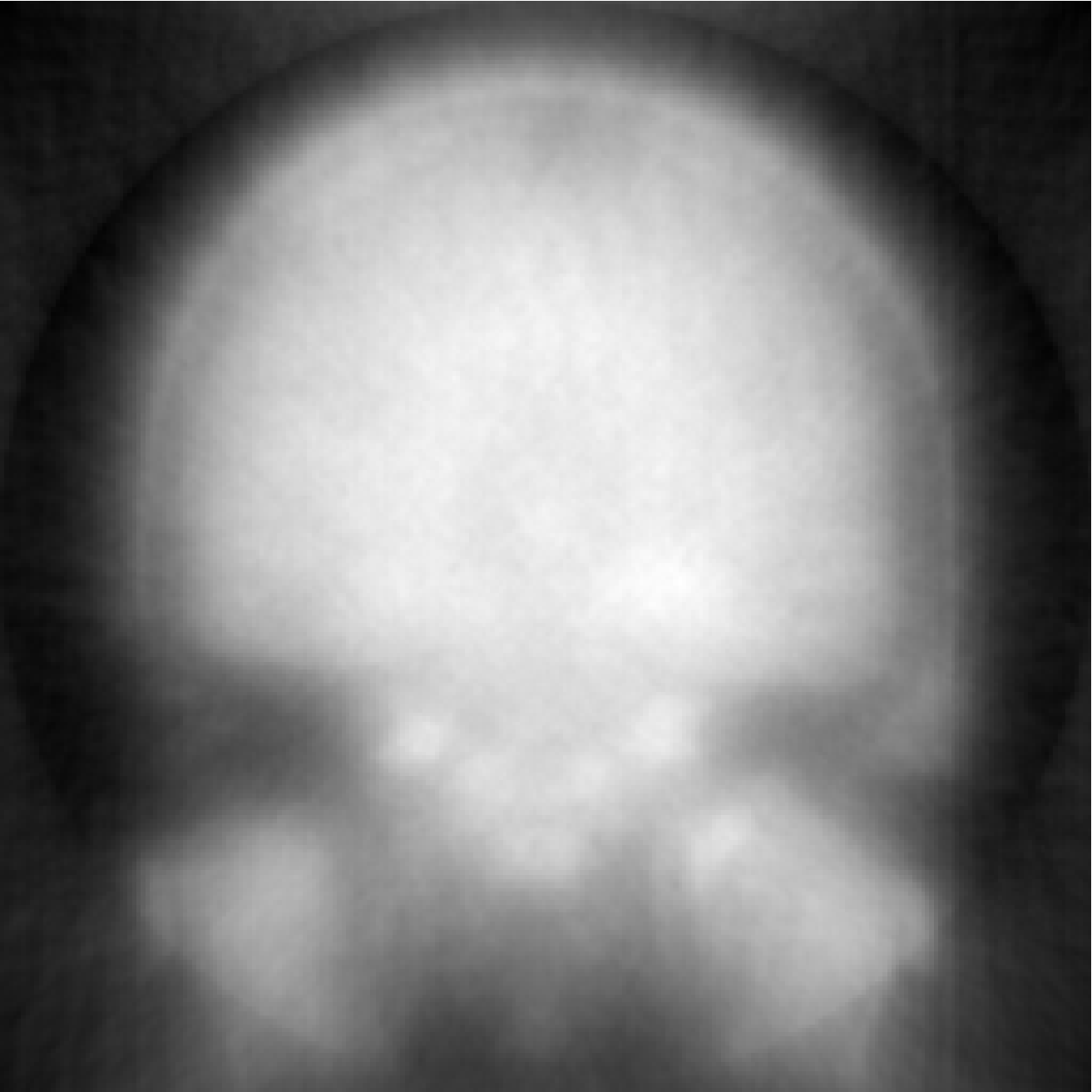}} &
 {\label{Fig:ABGKBGCVDP 18}\includegraphics[width=.20\textwidth]{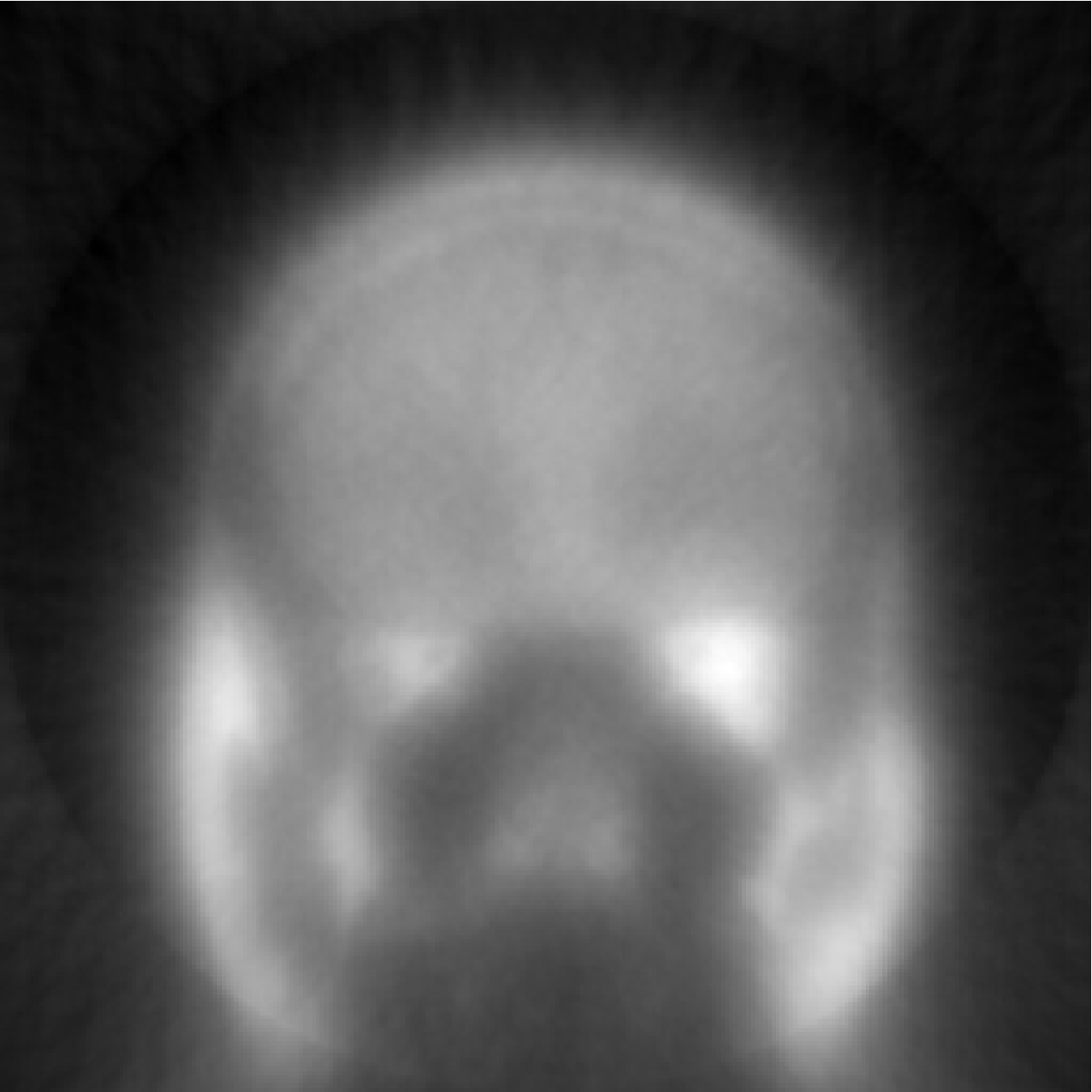}}  \\
  {\label{Fig:BAGKBGCVDP 5}\includegraphics[width=.20\textwidth]{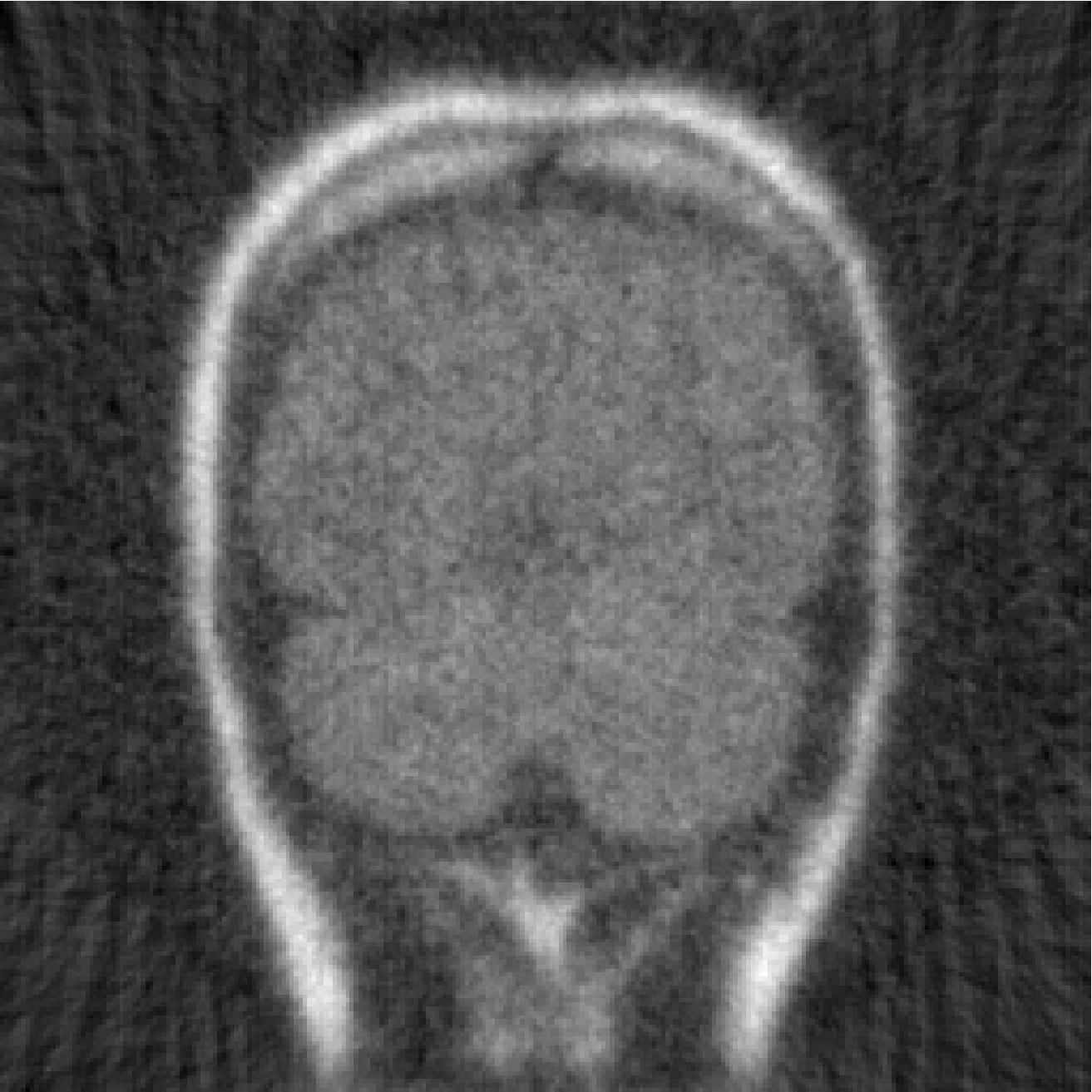}} &
 {\label{Fig:BAGKBGCVDP 11}\includegraphics[width=.20\textwidth]{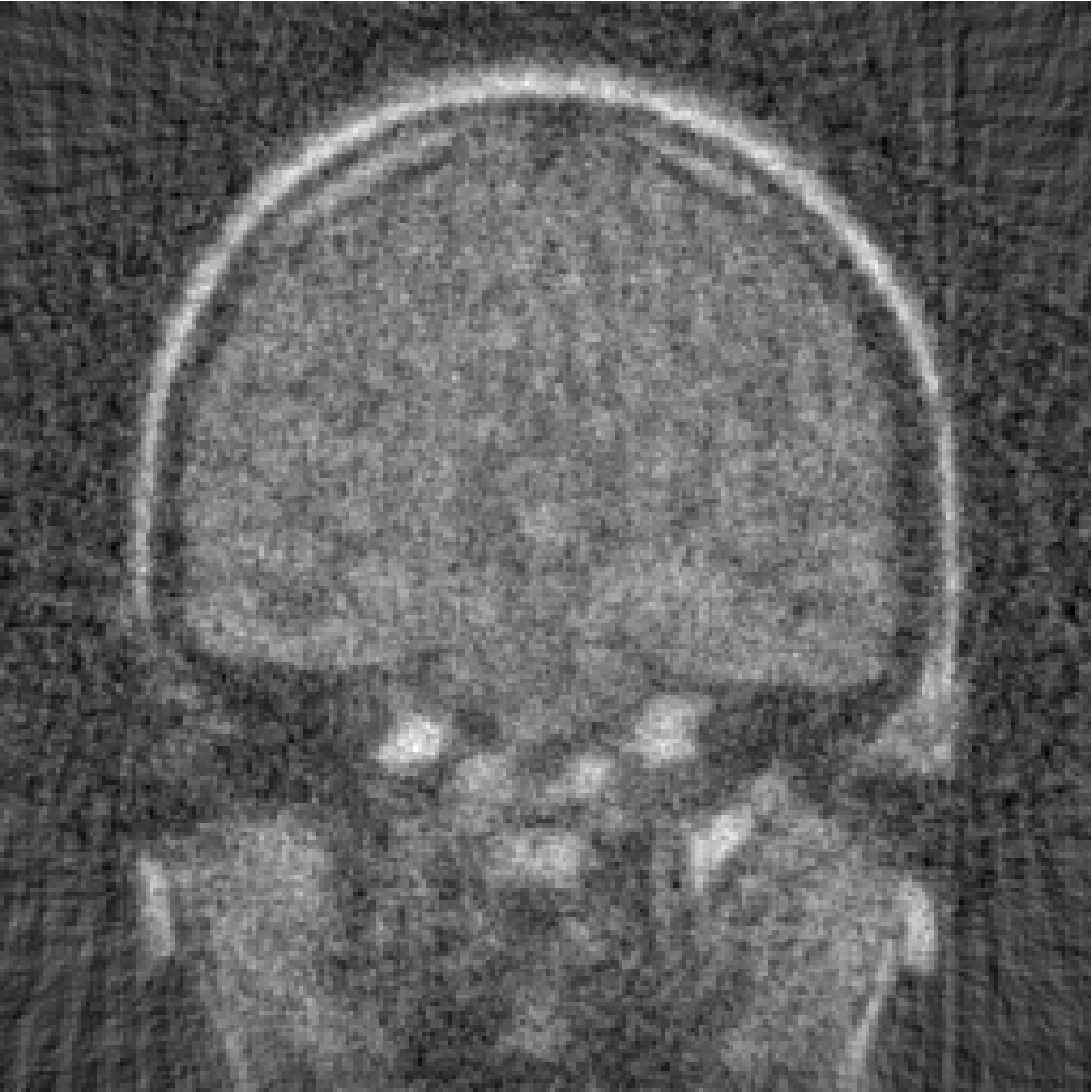}} &
 {\label{Fig:BAGKBGCVDP 18}\includegraphics[width=.20\textwidth]{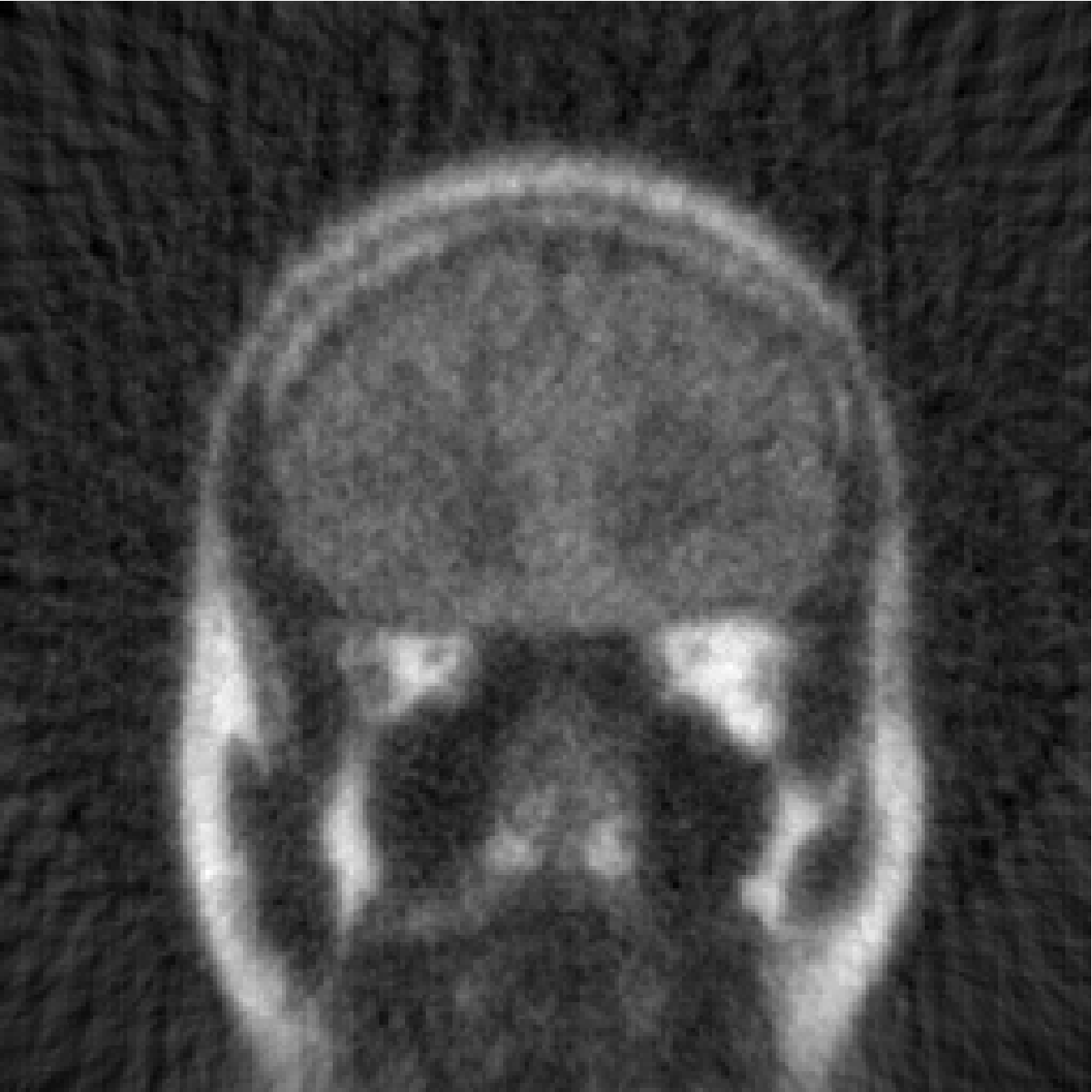}}  \\
 {\label{Fig:ABGMRESGCVDP 5}\includegraphics[width=.20\textwidth]{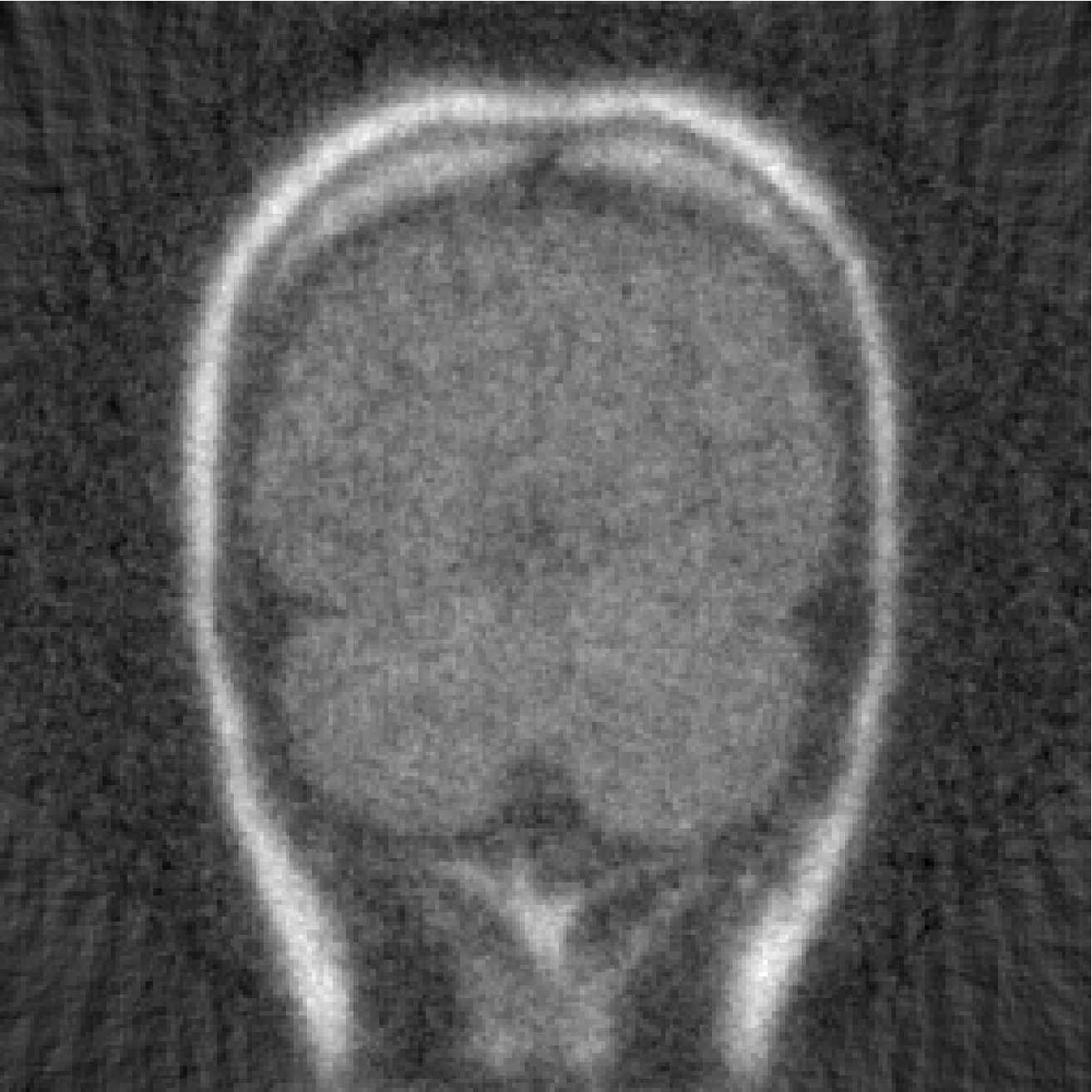}} &
 {\label{Fig:ABGMRESGCVDP 11}\includegraphics[width=.20\textwidth]{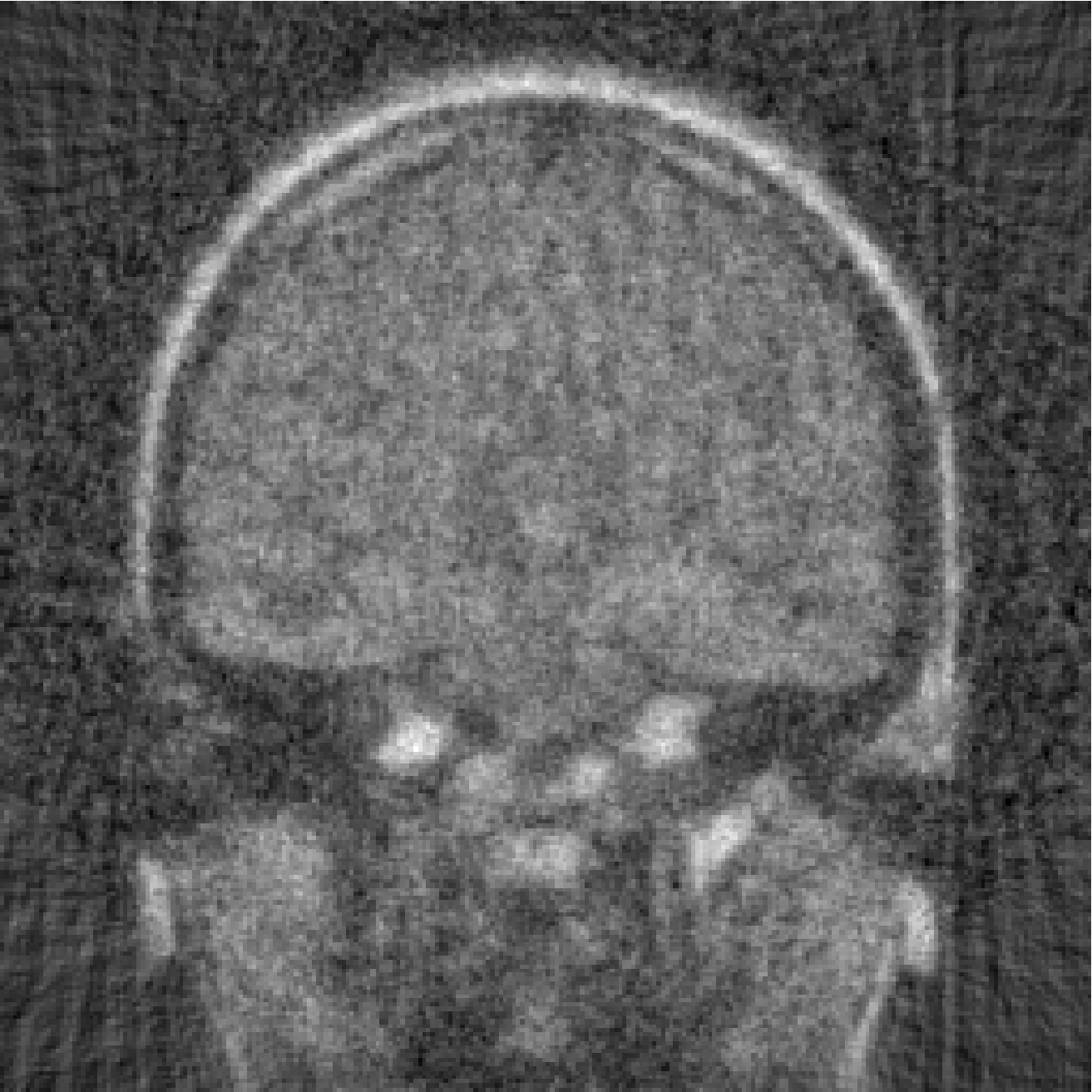}} &
 {\label{Fig:ABGMRESGCVDP 18}\includegraphics[width=.20\textwidth]{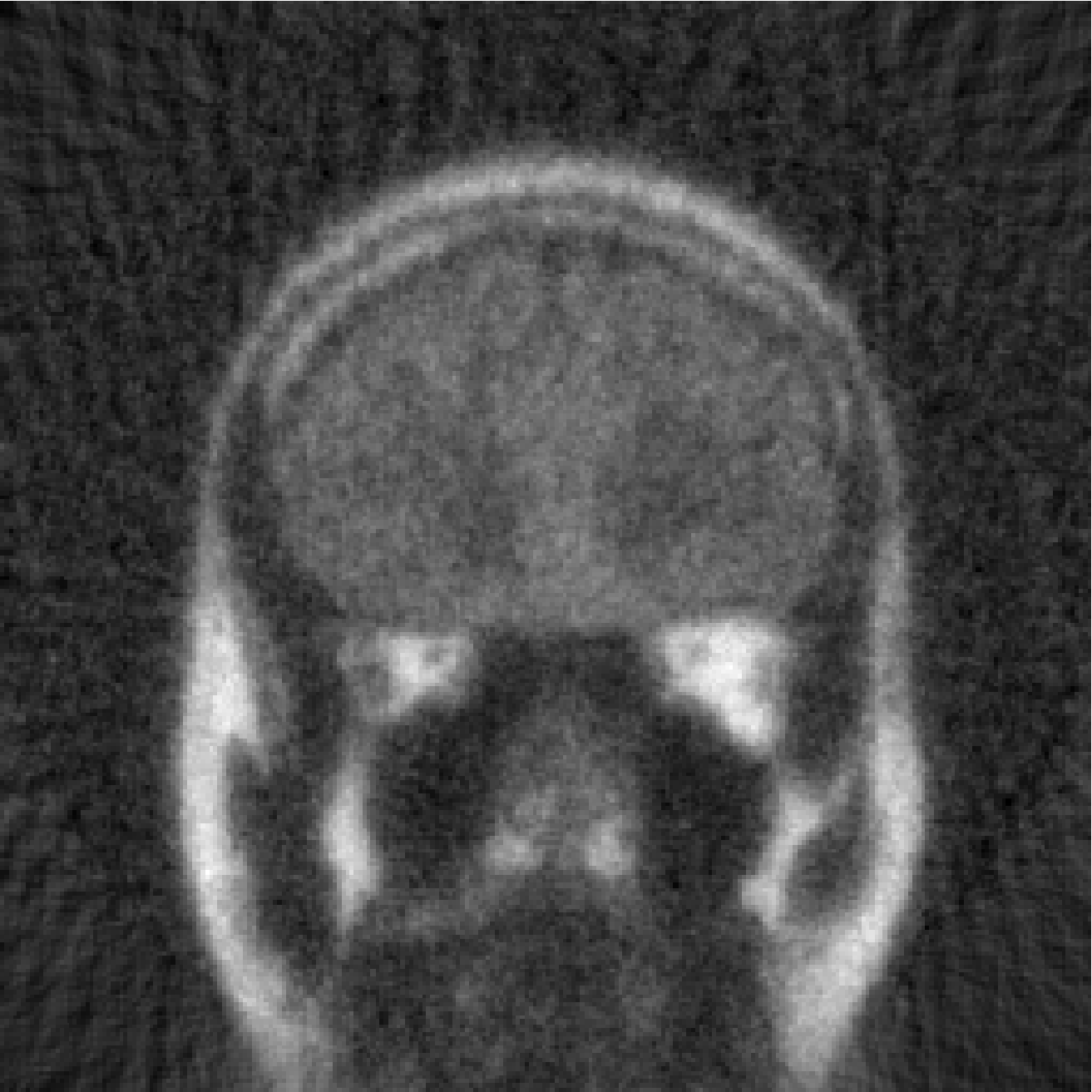}}  \\
  {\label{Fig:BAGMRESGCVDP 5}\includegraphics[width=.20\textwidth]{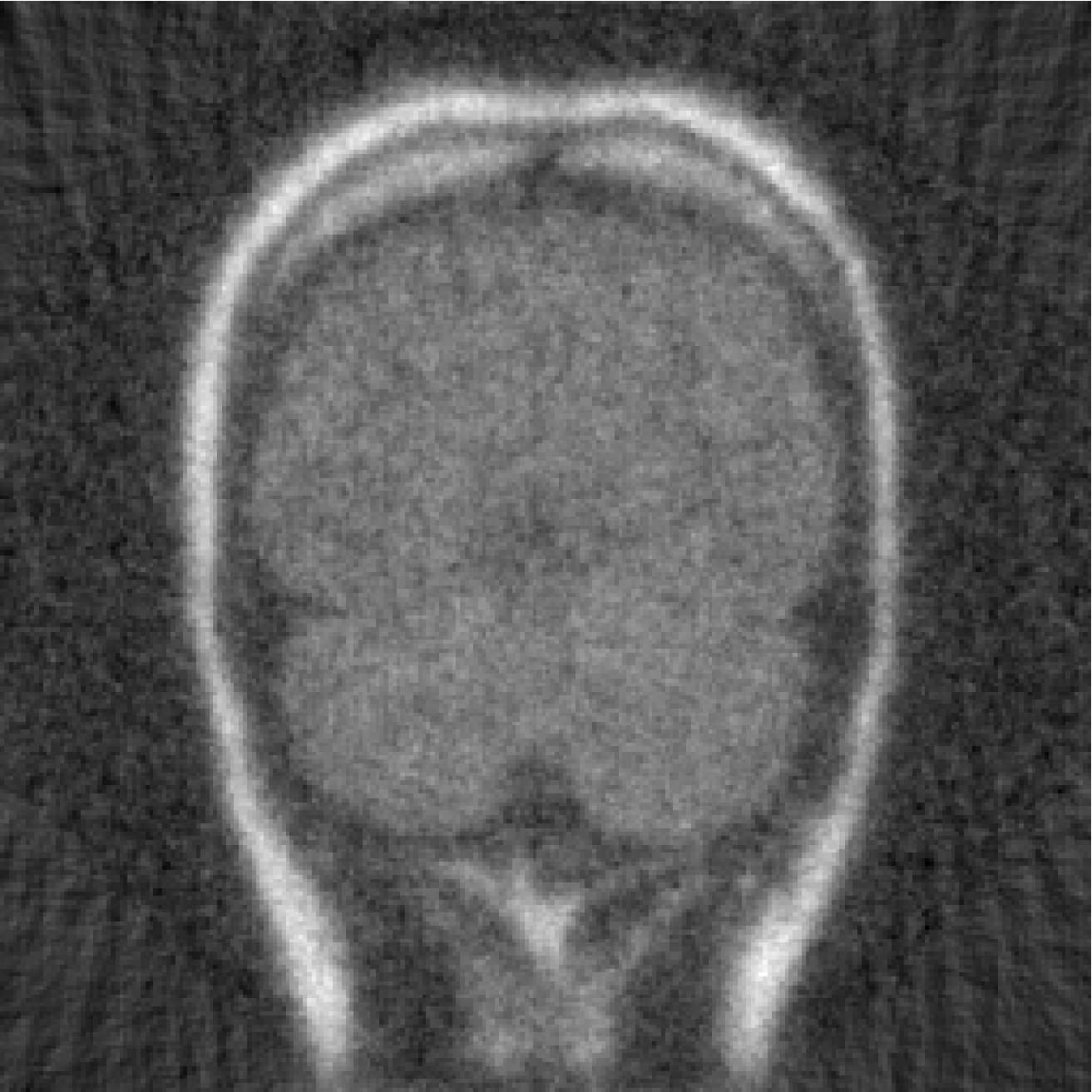}} &
 {\label{Fig:BAGMRESBGCVDP 11}\includegraphics[width=.20\textwidth]{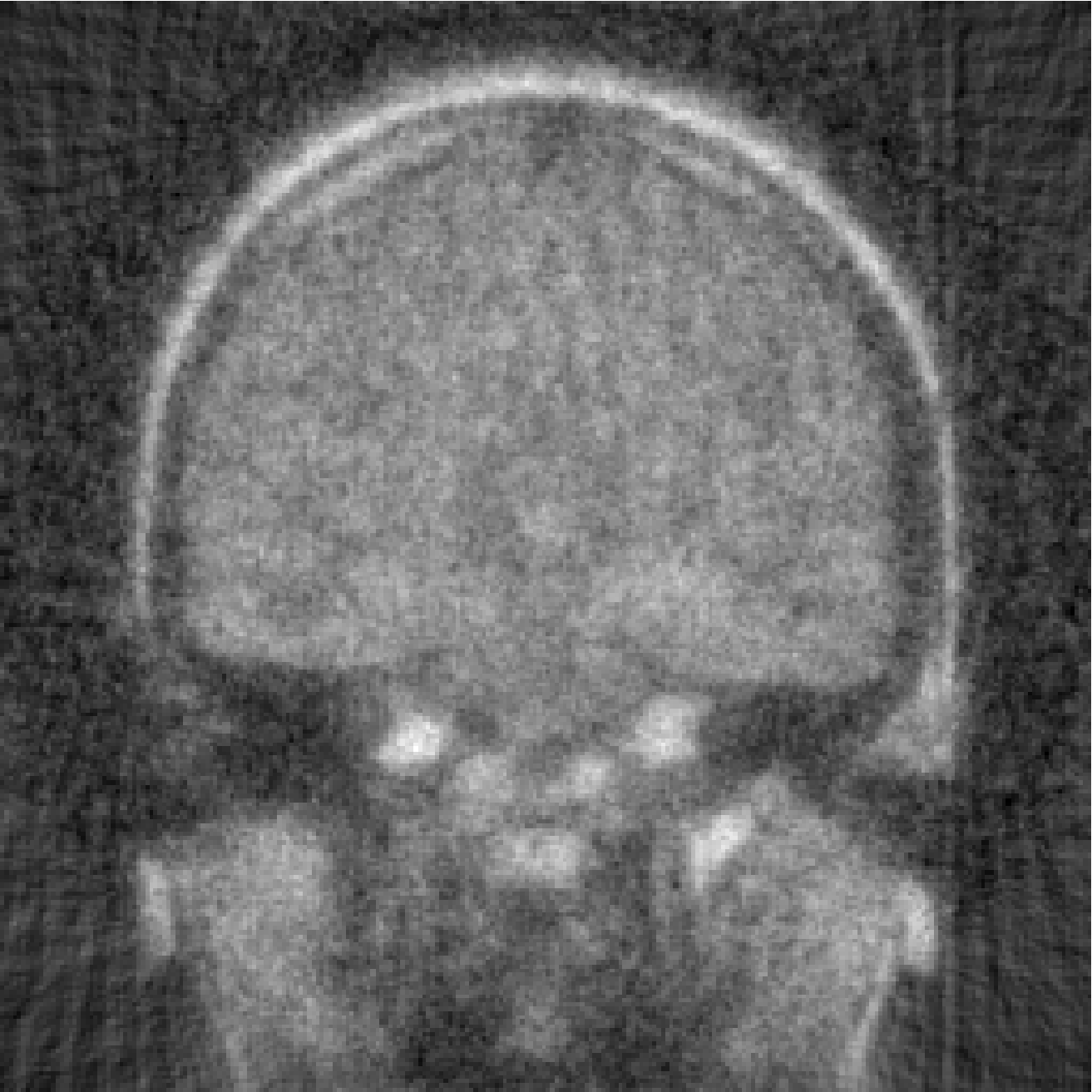}} &
 {\label{Fig:BAGMRESGCVDP 18}\includegraphics[width=.20\textwidth]{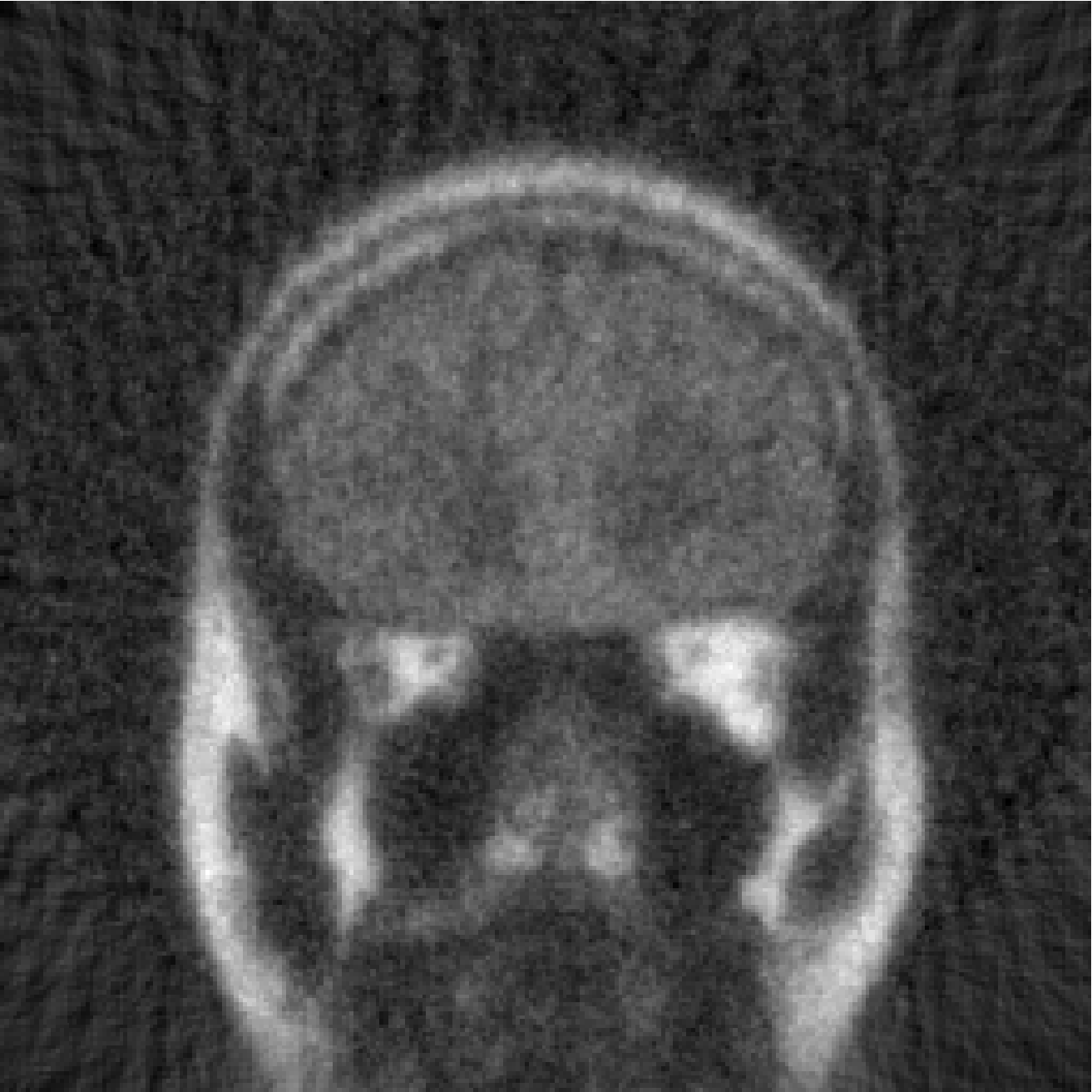}}  \\
  \end{tabular}
  \caption{Reconstructed images corresponding to \Cref{Fig:blurred slice 5,Fig:blurred slice 11,Fig:blurred slice 18}, obtained using hybrid AB-GKB, hybrid BA-GKB, hybrid AB-GMRES, and hybrid BA-GMRES with GCV and DP, shown in the first, second, third, and fourth rows, respectively.}\label{Fig:Recons slices}  
  \end{figure}

 \begin{example}\label{ex:problem2}[Under determined Example]
 We consider the MATLAB dataset $\texttt{mristack}$, which contains a 3D array of $21$ brain MRI slices, each of size $256\times 256$. Here, we consider the slices $5$, $11$, and $18$, as shown in \Cref{Fig:True slice 5,Fig:True slice 11,Fig:True slice 18}. Each slice is treated as a 2D X-ray tomography problem. We use the ASTRA Toolbox \citet{van2015astra} to generate the forward projector $A \in \mathbb{R}^{11520\times 65536}$  and backprojector $B \in \mathbb{R}^{65536\times 11520}$. The sinogram data are generated using parallel geometry through $45$ view angles. Then we use $8 \%$ white Gaussian noise, corresponding to SNR of approximately $21$, to obtain the contaminated sinogram data $\bfb \in \mathbb{R}^{11520}$, shown in \Cref{Fig:blurred slice 5,Fig:blurred slice 11,Fig:blurred slice 18}. The convergence results for slice $11$ are illustrated in \Cref{Fig:slices} for the GCV and L-curve, respectively, and quantitative results for all slices are reported in \Cref{Tab:slices Problem}. The reconstructed images using the GCV with the DP stopping rule for the different methods are shown in \Cref{Fig:Recons slices}. This example is designed to assess the performance of the proposed methods coupled with an underdetermined coefficient matrix.
 \end{example}

\begin{table}
 \caption{\label{Tab:slices Problem} The corresponding \text{RRE} \eqref{eq:RE} computed at $k_{DP}$ \eqref{eq:DP} and $k_{RNS}$ \eqref{eq:RNS} in parentheses, by all methods with $\lambda$ obtained with the GCV and L-curve, applied to the problem with $\text{SNR}\approx 21$ in \Cref{ex:problem2}. The lowest RRE results are shown in boldface. The timings are for all algorithms out to $150$ iterations and averaged over five runs.}
 \begin{tabular}{ccccc} 
 \toprule 
Method&Regularization&Iteration&$\text{RRE}(\bfx)$&Time(s)  \\ \midrule 
\multicolumn{5}{c}{Slice 5}\\     \hline 
\multirow{2}{*}{hybrid AB-GKB} & GCV& $150(3)$ &$0.56(0.65)$& \multirow{2}{*}{$1.84$}  \\ 
& L-curve & $72(3)$& $0.42(0.62)$ &  \\ 
\hline
\multirow{2}{*}{hybrid BA-GKB} & GCV& $19(3)$ &$\boldsymbol{0.33}(0.63)$& \multirow{2}{*}{$2.53$}  \\ 
& L-curve & $150(3)$& $0.48(0.63)$ &  \\ 
\hline
\multirow{2}{*}{hybrid AB-GMRES} & GCV& $6(11)$ &$\boldsymbol{0.33}(\boldsymbol{0.33})$& \multirow{2}{*}{$1.22$}  \\ 
& L-curve& $15(4)$& $0.41(0.59)$ &  \\ 
\hline
\multirow{2}{*}{hybrid BA-GMRES} & GCV& $8(8)$ &$\boldsymbol{0.33}(\boldsymbol{0.33})$& \multirow{2}{*}{$1.77$}  \\ 
& L-curve & $58(4)$& $0.50(0.52)$ & \\ 
\bottomrule
\multicolumn{5}{c}{Slice 11}\\     \bottomrule  
\multirow{2}{*}{hybrid AB-GKB} & GCV& $150(3)$ &$0.51(0.57)$& \multirow{2}{*}{$1.77$}  \\ 
& L-curve & $73(3)$& $0.47(0.56)$ &  \\ 
\hline
\multirow{2}{*}{hybrid BA-GKB} & GCV& $35(3)$ &$\boldsymbol{0.37}(0.57)$& \multirow{2}{*}{$2.41$}  \\ 
& L-curve & $150(3)$& $0.41(0.57)$ &  \\ 
\hline
\multirow{2}{*}{hybrid AB-GMRES} & GCV& $9(7)$ &$\boldsymbol{0.37}(0.54)$& \multirow{2}{*}{$1.16$}  \\ 
& L-curve& $16(3)$& $0.48(0.53)$ &  \\ 
\hline
\multirow{2}{*}{hybrid BA-GMRES} & GCV& $7(8)$ &$\boldsymbol{0.37}(\boldsymbol{0.36})$& \multirow{2}{*}{$1.76$}  \\ 
& L-curve & $66(5)$& $0.56(0.44)$ & \\ 
\bottomrule
\multicolumn{5}{c}{Slice 18}\\     \bottomrule  
\multirow{2}{*}{hybrid AB-GKB} & GCV& $150(3)$ &$0.61(0.72)$& \multirow{2}{*}{$1.73$}  \\ 
& L-curve & $70(6)$& $0.42(0.69)$ &  \\ 
\hline
\multirow{2}{*}{hybrid BA-GKB} & GCV& $21(7)$ &$\boldsymbol{0.35}(0.55)$& \multirow{2}{*}{$2.32$}  \\ 
& L-curve & $150(7)$& $0.57(0.55)$ &  \\ 
\hline
\multirow{2}{*}{hybrid AB-GMRES} & GCV& $7(15)$ &$\boldsymbol{0.35}(0.41)$& \multirow{2}{*}{$1.15$}  \\ 
& L-curve& $15(3)$& $0.41(0.69)$ &  \\ 
\hline
\multirow{2}{*}{hybrid BA-GMRES} & GCV& $9(10)$ &$\boldsymbol{0.35}(\boldsymbol{0.34})$& \multirow{2}{*}{$1.70$}  \\ 
& L-curve & $56(4)$& $0.49(0.60)$ & \\ \bottomrule 
 \end{tabular}    
 \end{table}

The results reported in \Cref{Tab:slices Problem} demonstrate the capabilities of the hybrid BA-GKB algorithms to yield results comparable to those obtained with hybrid GMRES methods. We observe that the AB-GKB reaches the maximum number of iterations without convergence when using the GCV but that the BA-GKB has the same issue when using the L-curve. The recorded time in column $5$ confirms the theoretical computational costs provided in \Cref{sec:comp costs}, which show that the GMRES methods are cheaper than the GKB methods. For this larger underdetermined ($m<<n$) problem, it is clear that the BA algorithms are slower than the AB algorithms, generally requiring more iterations to convergence and the cost being dominated by the $4nk^2>> 4mk^2$ term in \crefrange{cost:ABGMRES}{cost:BAGKB} in each case.

\section{Conclusions}\label{sec:Conclusion}
We introduced hybrid AB-GKB and BA-GKB algorithms for CT reconstruction with unmatched forward and backward projectors. The proposed methods were combined with GCV and L-curve parameter-selection strategies and evaluated using both the DP and RNS stopping criteria. Numerical results show that the proposed algorithms achieve reconstruction quality comparable to that of GMRES-based approaches. While GMRES remains computationally more efficient in terms of FLOPS, the proposed methods provide a viable alternative for solving unmatched-projector problems. No consistent superiority was observed between GCV and the L-curve methods, as their performance was problem-dependent. Moreover, the RNS stopping rule was generally satisfied faster than the DP stopping rule, but with a higher relative error.  Future work will investigate additional regularization strategies and their integration with the proposed hybrid frameworks for both GMRES and GKB formulations.



\bibliographystyle{cas-model2-names}

\bibliography{cas-refs}

@article{golub1965calculating,
  title={Calculating the singular values and pseudo-inverse of a matrix},
  author={Golub, Gene and Kahan, William},
  journal={Journal of the Society for Industrial and Applied Mathematics, Series B: Numerical Analysis},
  volume={2},
  number={2},
  pages={205--224},
  year={1965},
  publisher={SIAM}
}

@article{paige1982lsqr,
  title={{LSQR}: An algorithm for sparse linear equations and sparse least squares},
  author={Paige, Christopher C and Saunders, Michael A},
  journal={ACM Transactions on Mathematical Software (TOMS)},
  volume={8},
  number={1},
  pages={43--71},
  year={1982},
  publisher={ACM New York, NY, USA}
}

@article{ChungGazzolareview,
title = "Computational methods for large-scale inverse problems: a survey on hybrid projection methods",
author = "Julianne Chung and Silvia Gazzola",
year = "2024",
month = may,
day = "31",
volume = "66",
pages = "205--284",
journal = "SIAM Review",
number = "2",
}

@book{MATLAB,
year = {2026},
author = {MATLAB},
title = { Version 26.1.0.3251617 ({R}2026a)} }

@book{saad2003iterative,
  title={Iterative methods for sparse linear systems},
  author={Saad, Yousef},
  year={2003},
address={Philadelphia},
  publisher={SIAM}
}

@article{golub1979generalized,
  title={Generalized cross-validation as a method for choosing a good ridge parameter},
  author={Golub, Gene H and Heath, Michael and Wahba, Grace},
  journal={Technometrics},
  volume={21},
  number={2},
  pages={215--223},
  year={1979},
  publisher={Taylor \& Francis}
}

@article{hansen1993use,
  title={The use of the {L}-curve in the regularization of discrete ill-posed problems},
  author={Hansen, Per Christian and O’Leary, Dianne Prost},
  journal={SIAM Journal on Scientific Computing},
  volume={14},
  number={6},
  pages={1487--1503},
  year={1993},
  publisher={SIAM}
}

@article{hansen2022gmres,
  title={{GMRES} methods for tomographic reconstruction with an unmatched back projector},
  author={Hansen, Per Christian and Hayami, Ken and Morikuni, Keiichi},
  journal={Journal of Computational and Applied Mathematics},
  volume={413},
  pages={114352},
  year={2022},
  publisher={Elsevier}
}

@article{article,
author = {Zeng, Gengsheng and Gullberg, Grant},
year = {2000},
month = {06},
pages = {548-55},
title = {Unmatched Projector/Backprojector Pairs in an Iterative Reconstruction Algorithm},
volume = {19},
journal = {IEEE transactions on medical imaging},
doi = {10.1109/42.870265}
}

@article{doi:10.1137/070696313,
author = {Hayami, Ken and Yin, Jun-Feng and Ito, Tokushi},
title = {{GMRES} Methods for Least Squares Problems},
journal = {SIAM Journal on Matrix Analysis and Applications},
volume = {31},
number = {5},
pages = {2400-2430},
year = {2010},
doi = {10.1137/070696313}
}

@article{bentley2026hybrid,
  title={Hybrid {ABBA-GMRES} for Unmatched Backprojectors in Large Scale {X-Ray} Computerized Tomography},
  author={Bentley, Ryan and Pasha, Mirjeta and Landman, Malena Sabat{\'e} and Yang, Luisa and Zhang, Jeffery},
  journal={arXiv preprint arXiv:2602.17892},
  year={2026}
}

@article{doi:10.1137/17M1133828,
author = {Elfving, Tommy and Hansen, Per Christian},
title = {Unmatched Projector/Backprojector Pairs: Perturbation and Convergence Analysis},
journal = {SIAM Journal on Scientific Computing},
volume = {40},
number = {1},
pages = {A573-A591},
year = {2018},
doi = {10.1137/17M1133828}
}

@inproceedings{hansen2021stopping,
  title={Stopping rules for algebraic iterative reconstruction methods in computed tomography},
  author={Hansen, Per Christian and J{\o}rgensen, Jakob Sauer and Rasmussen, Peter Winkel},
  booktitle={2021 21st International Conference on Computational Science and its Applications (ICCSA)},
  pages={60--70},
  year={2021},
  organization={IEEE}
}

@book{morozov2012methods,
  title={Methods for solving incorrectly posed problems},
  author={Morozov, Vladimir Alekseevich},
  year={2012},
  publisher={Springer Science \& Business Media},
  address={New York}
}

@article{van2015astra,
  title={{The ASTRA Toolbox: A platform for advanced algorithm development in electron tomography}},
  author={Van Aarle, Wim and Palenstijn, Willem Jan and De Beenhouwer, Jan and Altantzis, Thomas and Bals, Sara and Batenburg, K Joost and Sijbers, Jan},
  journal={Ultramicroscopy},
  volume={157},
  pages={35--47},
  year={2015},
  publisher={Elsevier}
}

@article{fong2011lsmr,
  title={{LSMR}: An iterative algorithm for sparse least-squares problems},
  author={Fong, David Chin-Lung and Saunders, Michael},
  journal={SIAM Journal on Scientific Computing},
  volume={33},
  number={5},
  pages={2950--2971},
  year={2011},
  publisher={SIAM}
}

@article{bianchi2026iterated,
  title={The iterated {Golub-Kahan-Tikhonov} method},
  author={Bianchi, Davide and Donatelli, Marco and Furch{\'\i}, Davide and Reichel, Lothar},
  journal={BIT Numerical Mathematics},
  volume={66},
  number={1},
  pages={18},
  year={2026},
  publisher={Springer}
}

@article{beik2020golub,
  title={{Golub--Kahan bidiagonalization for ill-conditioned tensor equations with applications}},
  author={Beik, Fatemeh PA and Jbilou, Khalide and Najafi-Kalyani, Mehdi and Reichel, Lothar},
  journal={Numerical Algorithms},
  volume={84},
  number={4},
  pages={1535--1563},
  year={2020},
  publisher={Springer}
}

@article{bentbib2018solution,
  title={Solution methods for linear discrete ill-posed problems for color image restoration},
  author={Bentbib, Abdeslem H and El Guide, M and Jbilou, Khalide and Onunwor, E and Reichel, Lothar},
  journal={BIT Numerical Mathematics},
  volume={58},
  number={3},
  pages={555--576},
  year={2018},
  publisher={Springer}
}

@article{bjorck1988bidiagonalization,
  title={A bidiagonalization algorithm for solving large and sparse ill-posed systems of linear equations},
  author={Bj{\"o}rck, {\AA}ke},
  journal={BIT Numerical Mathematics},
  volume={28},
  number={3},
  pages={659--670},
  year={1988},
  publisher={Springer}
}

@article{alsubhi2026GKB,
  title={{GKB Methods for X-Ray Computed Tomography with an Unmatched Back Projector}},
  author={Alsubhi, Abdulmajeed },
  journal={arXiv preprint arxiv.org/abs/2606.31153},
  year={2026}
}

@article{varah1983pitfalls,
  title={Pitfalls in the numerical solution of linear ill-posed problems},
  author={Varah, James M},
  journal={SIAM Journal on Scientific and Statistical Computing},
  volume={4},
  number={2},
  pages={164--176},
  year={1983},
  publisher={SIAM}
}

@article{wathen2025least,
  title={Least Squares and the Not-Normal Equations},
  author={Wathen, Andrew J},
  journal={SIAM Review},
  volume={67},
  number={4},
  pages={865--872},
  year={2025},
  publisher={SIAM}
}






\end{document}